\documentclass[a4paper,12pt]{amsart}
\usepackage[top=3cm,left=2.5cm,right=2.5cm,bottom=3cm]{geometry}
\usepackage{relsize}
\usepackage{lineno}
\usepackage{hyperref}
\usepackage{amsmath,amsthm,pb-diagram,amssymb,comment}
\usepackage{amsfonts,graphicx,color}
\usepackage{enumerate, fancyhdr, dsfont}
\usepackage{thmtools,hyperref,cleveref}
\usepackage[normalem]{ulem}
\usepackage{stmaryrd}
\usepackage{textcomp}
\usepackage{mathabx}
\usepackage{tikz,float}
\hypersetup{
    colorlinks=true,
    linkcolor=dodger,
    filecolor=dodger,
    urlcolor=dodger,
    citecolor=dodger,
}

\usepackage{caption}
\usepackage{booktabs}
\usepackage{mwe}

\newcommand*\DrawHoriLine{%
  \par
  \dimen@ \prevdepth
  \kern -.1\p@
  \hrule \@height .1\p@ \@depth .1\p@
  \kern -.1\p@
  \prevdepth \z@
  \null
  \vskip -\baselineskip
  \prevdepth \dimen@
}

    \newcommand{\ran}{\mbox{\rm ran}}

    \newcommand{\thzfc}{\mathrm{ZFC}}

     \newcommand{\Ed}{\mathbf{Ed}}

    \newcommand{\Cn}{\mathbf{Cn}}

    \newcommand{\Cwf}{\mathcal{C}}
    
    \newcommand{\Ewf}{\mathcal{E}}
    \newcommand{\Fwf}{\mathcal{F}}
    
    \newcommand{\Hwf}{\mathcal{H}}
    \newcommand{\Iwf}{\mathcal{I}}
    \newcommand{\Jwf}{\mathcal{J}}
    \newcommand{\Mwf}{\mathcal{M}}
    \newcommand{\Nwf}{\mathcal{N}}
    
    \newcommand{\Pwf}{\mathcal{P}}
    \newcommand{\Swf}{\mathcal{S}}

    \newcommand{\bfrak}{\mathfrak{b}}
    \newcommand{\cfrak}{\mathfrak{c}}
    \newcommand{\dfrak}{\mathfrak{d}}

    \newcommand{\ufrak}{\mathfrak{u}}

    \newcommand{\menos}{\smallsetminus}

    \newcommand{\frestr}{\!\!\upharpoonright\!\!}

    \newcommand{\add}{\mbox{\rm add}}
    \newcommand{\cov}{\mbox{\rm cov}}
    \newcommand{\non}{\mbox{\rm non}}
    \newcommand{\cof}{\mbox{\rm cof}}
    \newcommand{\limdir}{\mbox{\rm limdir}}

    \newcommand{\Bor}{\mathds{B}}
    \newcommand{\Cor}{\mathds{C}}
    \newcommand{\Dor}{\mathds{D}}
    \newcommand{\Eor}{\mathds{E}}

    \newcommand{\LOCor}{\mathds{LOC}}

    \newcommand{\Por}{\mathds{P}}
    \newcommand{\Pbb}{\mathds{P}}
    \newcommand{\Qor}{\mathds{Q}}
    \newcommand{\Ior}{\mathds{I}}

    \newcommand{\Qnm}{\dot{\mathds{Q}}}

    \newcommand{\SNwf}{\mathcal{SN}}

    \newcommand{\R}{\mathbb{R}}

    \newcommand{\cf}{\mbox{\rm cf}}

    \newcommand{\la}{\langle}
    \newcommand{\ra}{\rangle}

\newcommand{\Dbf}{\mathbf{D}}
\newcommand{\mbf}{\mathbf{m}}

\newcommand{\uf}{\mathrm{uf}}
\newcommand{\seq}{\mathrm{seq}}

\newcommand{\Fr}{\mathrm{Fr}}
\newcommand{\Rbf}{\mathbf{R}}

\newcommand{\Cbf}{\mathbf{C}}
\newcommand{\Lc}{\mathbf{Lc}}

\newcommand{\loc}{\mathbf{Lc}}

\newcommand{\Lb}{\mathbf{Lb}}
\newcommand{\Hcal}{\mathcal{H}}
\newcommand{\Scal}{\mathcal{S}}
\newcommand{\id}{\mathrm{id}}
\newcommand{\blc}{\mathfrak{b}^{\mathrm{Lc}}}
\newcommand{\dlc}{\mathfrak{d}^{\mathrm{Lc}}}
\newcommand{\dloc}{\mathfrak{d}^{\mathrm{Lc_*}}}
\newcommand{\bloc}{\mathfrak{b}^{\mathrm{Lc_*}}}
\newcommand{\balc}{\mathfrak{b}^{\mathrm{aLc}}}
\newcommand{\dalc}{\mathfrak{d}^{\mathrm{aLc}}}

\newcommand{\vfa}{\mathfrak{b}^{\mathrm{Lc}}}
\newcommand{\cfa}{\mathfrak{d}^{\mathrm{Lc}}}

\newcommand{\leqT}{\preceq_{\mathrm{T}}}
\newcommand{\eqT}{\cong_{\mathrm{T}}}

\newcommand{\Ue}{\mathbf{Ue}}

\newcommand{\onemainP}{\textrm{(A1)}}
\newcommand{\twomainP}{\textrm{(A2)}}

\newcommand{\mn}{{\medskip\noindent}}
\newcommand{\sn}{{\smallskip\noindent}}

\newenvironment{PROOF}[2][\proofname.]
   {\begin{proof}[#1]\renewcommand{\qedsymbol}{\textsquare$_{\rm #2}$}}
   {\end{proof}}

\DeclareSymbolFont{extraup}{U}{zavm}{m}{n}
\DeclareMathSymbol{\varheart}{\mathalpha}{extraup}{86}
\DeclareMathSymbol{\vardiamond}{\mathalpha}{extraup}{87}

\definecolor{dodger}{rgb}{0.0,0.5,1.0}

\definecolor{amber}{rgb}{1.0,0.49,0.0}

\definecolor{ogreen}{RGB}{107,142,35}

\title[Closed measure zero sets]{A friendly iteration forcing that the four cardinal characteristics of $\Ewf$ can be pairwise different}

\author{Miguel A. Cardona}
\address{Pavol Jozef \v{S}af\'arik University}
\email{miguel.cardona@upjs.sk}
\urladdr{https://www.researchgate.net/profile/Miguel\_Cardona\_Montoya}



\begin{document}

\makeatletter
\def\@roman#1{\romannumeral #1}
\makeatother

\newcounter{enuAlph}
\renewcommand{\theenuAlph}{\Alph{enuAlph}}

\renewcommand{\theequation}{\thesection.\arabic{equation}}

\theoremstyle{plain}
  \newtheorem{theorem}{Theorem}[section]
  \newtheorem{corollary}[theorem]{Corollary}
  \newtheorem{lemma}[theorem]{Lemma}
  \newtheorem{mainlemma}[theorem]{Main Lemma}
  \newtheorem{prop}[theorem]{Proposition}
  \newtheorem{clm}[theorem]{Claim}
  \newtheorem{fact}[theorem]{Fact}
  \newtheorem{exer}[theorem]{Exercise}
  \newtheorem{question}[theorem]{Question}
  \newtheorem{problem}[theorem]{Problem}
  \newtheorem{conjecture}[theorem]{Conjecture}
  \newtheorem*{thm}{Theorem}
  \newtheorem*{maintheorem}{Main Theorem}
  \newtheorem{mainproblem}[enuAlph]{Main Problem}
  \newtheorem{question*}[enuAlph]{Question}
  \newtheorem{teorema}[enuAlph]{Theorem}
  \newtheorem*{corolario}{Corollary}
\theoremstyle{definition}
  \newtheorem{definition}[theorem]{Definition}
  \newtheorem{example}[theorem]{Example}
  \newtheorem{remark}[theorem]{Remark}
  \newtheorem{notation}[theorem]{Notation}
  \newtheorem{context}[theorem]{Context}

  \newtheorem*{defi}{Definition}
  \newtheorem*{acknowledgements}{Acknowledgements}

  \def\sectionautorefname{Section}
\def\subsectionautorefname{Subsection}

\maketitle


\begin{abstract}
Let $\Ewf$ be the $\sigma$-ideal generated by the closed measure zero sets of reals. We use an ultrafilter-extendable matrix iteration of ccc posets to force that, for $\Ewf$, their associated cardinal characteristics (i.e.\ additivity, covering, uniformity and cofinality) are pairwise different. 
\end{abstract}

\makeatother


\section{Introduction}


Let $\Mwf$ and $\Nwf$, as usual, denote the $\sigma$-ideal of first category subsets of $\R$ and the $\sigma$-ideal of Lebesgue null subsets of $\R$ respectively and let $\Ewf$ be the $\sigma$-ideal generated by closed measure zero subsets of $\R$. It is well-known that $\Ewf\subseteq\Nwf\cap\Mwf$. Even more, it was proved that $\Ewf$ is a proper subideal of $\Nwf\cap\Mwf$ (see~\cite[Lemma 2.6.1]{BJ}).

This work belongs to the framework of consistency results where several cardinal characteristics are pairwise different. In particular,  we focus on the four cardinal characteristics  associated with $\Ewf$. In general, those cardinals are defined as follows: Let $\Iwf$ be an ideal of subsets of $X$ such that $\{x\}\in \Iwf$ for all $x\in X$. We define \emph{the cardinal characteristics associated with $\Iwf$} by

\noindent\textbf{Additivity of $\Iwf$:} $\add(\Iwf)=\min\{|\Jwf|:\,\Jwf\subseteq\Iwf,\,\bigcup\Jwf\notin\Iwf\}$. 
\smallskip

\noindent\textbf{Covering of $\Iwf$:} $\cov(\Iwf)=\min\{|\Jwf|:\,\Jwf\subseteq\Iwf,\,\bigcup\Jwf=X\}$.
\smallskip

\noindent\textbf{Uniformity of $\Iwf$:} $\non(\Iwf)=\min\{|A|:\,A\subseteq X,\,A\notin\Iwf\}$.
\smallskip

\noindent\textbf{Cofinality of $\Iwf$:} $\cof(\Iwf)=\min\{|\Jwf|:\,\Jwf\subseteq\Iwf,\ \forall\, A\in\Iwf\ \exists\, B\in \Jwf\colon A\subseteq B\}$.

\begin{figure}[h]
  \begin{center}
    \includegraphics[scale=1.0]{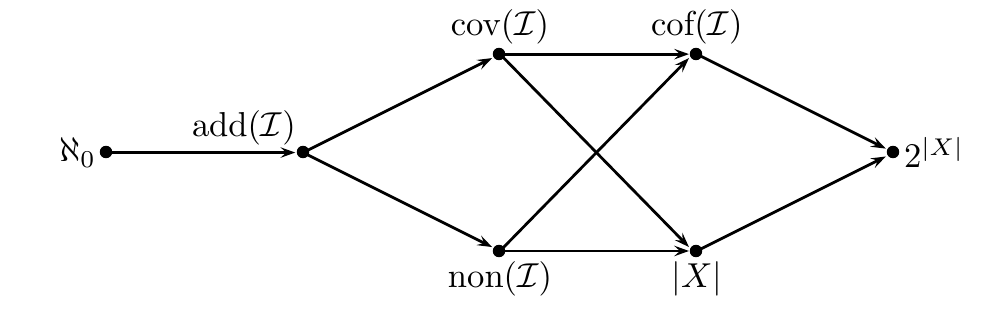}
    \caption{Diagram of the cardinal characteristics associated with $\Iwf$. An arrow  $\mathfrak x\rightarrow\mathfrak y$ means that (provably in ZFC) 
    $\mathfrak x\le\mathfrak y$.}
    \label{diag:idealI}
  \end{center}
\end{figure}

For ideals $\Iwf\subseteq\Jwf$ define 
\[\add(\Iwf, \Jwf)=\min\{|\Fwf|:\Fwf\subseteq\Iwf\textrm{\ and\ } \bigcup\Fwf\notin\Jwf\}\]
and 
\[\cof(\Iwf, \Jwf)=\min\{|\Fwf|:\!\Fwf\subseteq\Jwf\textrm{\ and\ }(\forall A\in\Iwf)(\exists B\in \Fwf)(A\subseteq B)\}.\]

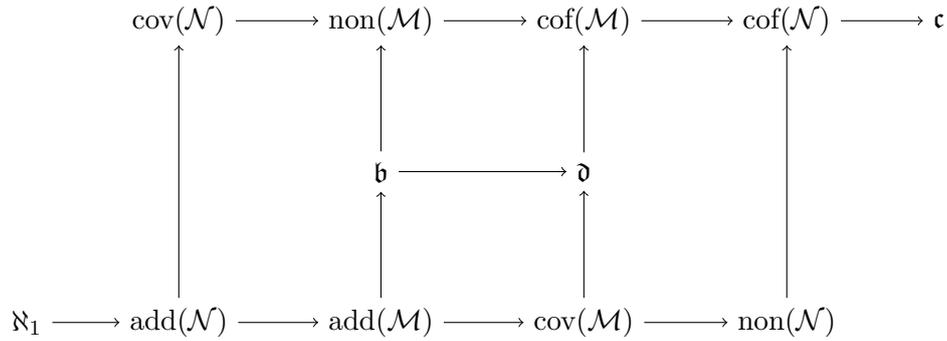
\begin{figure}[ht!]
\centering
\begin{tikzpicture}
\small{
\node (aleph1) at (-1,3) {$\aleph_1$};
\node (addn) at (1,3){$\add(\Nwf)$};
\node (covn) at (1,7){$\cov(\Nwf)$};
\node (nonn) at (9,3) {$\non(\Nwf)$} ;
\node (cfn) at (9,7) {$\cof(\Nwf)$} ;
\node (addm) at (3.66,3) {$\add(\Mwf)$} ;
\node (covm) at (6.33,3) {$\cov(\Mwf)$} ;
\node (nonm) at (3.66,7) {$\non(\Mwf)$} ;
\node (cfm) at (6.33,7) {$\cof(\Mwf)$} ;
\node (b) at (3.66,5) {$\bfrak$};
\node (d) at (6.33,5) {$\dfrak$};
\node (c) at (11,7) {$\cfrak$};
\draw (aleph1) edge[->] (addn)
      (addn) edge[->] (covn)
      (covn) edge [->] (nonm)
      (nonm)edge [->] (cfm)
      (cfm)edge [->] (cfn)
      (cfn) edge[->] (c);
\draw
   (addn) edge [->]  (addm)
   (addm) edge [->]  (covm)
   (covm) edge [->]  (nonn)
   (nonn) edge [->]  (cfn);
\draw (addm) edge [->] (b)
      (b)  edge [->] (nonm);
\draw (covm) edge [->] (d)
      (d)  edge[->] (cfm);
\draw (b) edge [->] (d);


}
\end{tikzpicture}
\caption{Cicho\'n's diagram. The arrows mean $\leq$ and dotted arrows represent
  $\add(\Mwf)=\min\{\bfrak,\cov(\Mwf)\}$ and $\cof(\Mwf)=\max\{\dfrak,\non(\Mwf)\}$.}
  \label{FigCichon}
\end{figure}

\autoref{diag:idealI} shows the natural inequalities between the cardinal characteristics associated with  $\Iwf$. They have been studied intensively for the ideal $\Ewf$. In particular, they are closely related to the cardinal characteristics in Cicho\'n's diagram (see~\autoref{FigCichon}). Some of the earliest results are due to Miller \cite{Miller}, who proved $\add(\Ewf, \Mwf)\leq\dfrak$ and $\cof(\Ewf, \Mwf)\geq\dfrak$. 

Here, $\bfrak$ and $\dfrak$ denotes, as usual, the classical \emph{unbounded number} and \emph{dominating number}, respectively (see~\autoref{def:b_d}).

On the other hand, Bartoszy\'nski and Shelah~\cite{BS1992} proved the following in ZFC, which is probably the deepest result about the cardinal characteristics of $\Ewf$.

\begin{theorem}[{\cite{BS1992}, see also \cite[Thm.~2.6.17]{BJ}}]\label{bartoS}
$\add(\Ewf)=\add(\Mwf)$ and $\cof(\Ewf)=\cof(\Mwf)$.
\end{theorem}


\begin{figure}[!htb]
\includegraphics[scale=1]{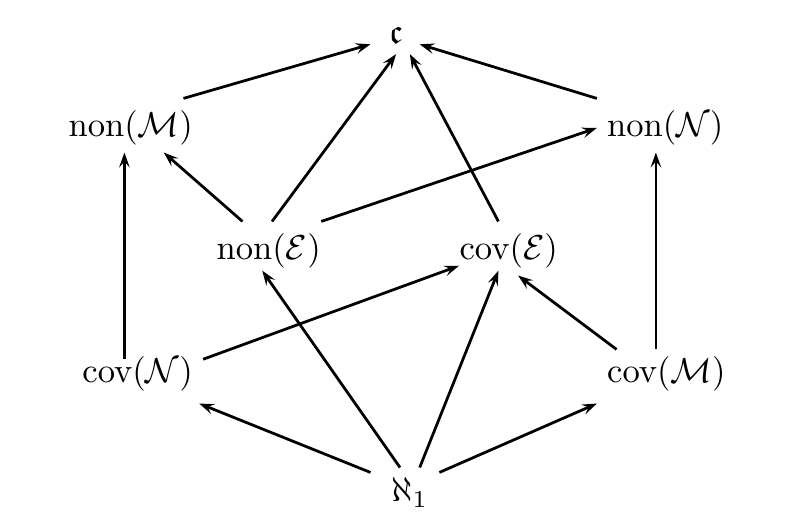}
\caption{The arrows mean that $\leq$ is provable in ZFC.}
\label{diagE}
\end{figure}

Furthermore, using the fact that $\Ewf$ is contained in the intersection of $\Mwf$ and $\Nwf$, the inequalities $\cov(\Ewf)\geq\max\{\cov(\Mwf),\cov(\Nwf)\}$ and $\non(\Ewf)\leq\min\{\non(\Mwf),\non(\Nwf)\}$ are easy to
establish (see~\autoref{diagE}). We denote, as usual, $\cfrak:=2^{\aleph_0}=|\R|$, and recall that $\aleph_1$ is the smallest uncountable cardinal.

The answer to the following is not known. 

\begin{question*}\label{conjE}
Is it consistent that all the cardinals in~\autoref{diagE} are pairwise different?
\end{question*}

Bartoszy\'nski and Shelah~\cite[Sec.~5]{BS1992} obtained a model by adding $\aleph_2$-many random reals by using a finite support (FS) iteration of lenght $\aleph_2$ over a model of CH that satisfies $\non(\Ewf)=\aleph_1<\cov(\Nwf)=\non(\Mwf)=\cov(\Mwf)=\non(\Nwf)=\cov(\Ewf)=\aleph_2$. We point out the following easy consequence which involve the latter result due mostly to~Bartoszy\'nski and Shelah~\cite{BS1992}:

\begin{teorema}{\emph{(\autoref{weakversionE})}}\label{WVcovnonE}
Let $\theta\leq \nu$ be uncountable regular cardinals and let $\lambda$ be a cardinal such that $\nu\leq\lambda=\lambda^{{<}\theta}$. Then there a ccc poset forcing
\[\non(\Ewf)=\bfrak=\theta\leq\cov(\Nwf)=\non(\Mwf)=\cov(\Mwf)=\non(\Nwf)=\nu\leq\cov(\Ewf)=\dfrak=\lambda.\]
\end{teorema}

Regarding~\autoref{diagE} again, we present a theorem of Bartoszy\'nski and Shelah which gives an upper bound for $\cov(\Ewf)$ and an lower bound for $\non(\Ewf)$ different to $\cfrak$ and $\aleph_1$, respectively.

\begin{theorem}[{\cite{BS1992}, see also \cite[Thm.~2.6.9]{BJ}}]\label{covnonE} \textrm{}
\begin{itemize}
    \item[(1)]  $\max\{\cov(\Mwf),\cov(\Nwf)\}\leq\cov(\Ewf)\leq\max\{\dfrak,\cov(\Nwf)\}$.
    \mn
    \item[(2)]   $\min\{\bfrak,\non(\Nwf)\}\leq\non(\Ewf)\leq\min\{\non(\Mwf),\non(\Nwf)\}$. 
\end{itemize}
In particular, 
\begin{itemize}
    \item[(i)]  if $\dfrak=\cov(\Mwf)$, then $\cov(\Ewf)=\max\{\cov(\Mwf),\cov(\Nwf)\}$,\mn
    \item[(ii)]  if $\bfrak=\non(\Mwf)$, then $\non(\Ewf)=\min\{\non(\Mwf),\non(\Nwf)\}$. 
\end{itemize}
\end{theorem}

Concerning~\autoref{FigCichon}, a natural question to ask is the following.

\begin{question*}\label{Conj:Cichons}
Is it consistent that all the cardinals in~\autoref{FigCichon} (with the exception of the dependent values $\add(\Mwf)$ and $\cof(\Mwf)$) are pairwise different?
\end{question*}

\begin{figure}[ht!]
\centering
\begin{tikzpicture}
\small{
\node (aleph1) at (-1,3) {$\aleph_1$};
\node (addn) at (1,3){$\add(\Nwf)$};
\node (covn) at (1,7){$\cov(\Nwf)$};
\node (nonn) at (9,3) {$\non(\Nwf)$} ;
\node (cfn) at (9,7) {$\cof(\Nwf)$} ;
\node (addm) at (3.66,3) {$\add(\Mwf)$} ;
\node (covm) at (6.33,3) {$\cov(\Mwf)$} ;
\node (nonm) at (3.66,7) {$\non(\Mwf)$} ;
\node (cfm) at (6.33,7) {$\cof(\Mwf)$} ;
\node (b) at (3.66,5) {$\bfrak$};
\node (d) at (6.33,5) {$\dfrak$};
\node (c) at (11,7) {$\cfrak$};
\draw (aleph1) edge[->] (addn)
      (addn) edge[->] (covn)
      (covn) edge [->] (nonm)
      (nonm)edge [->] (cfm)
      (cfm)edge [->] (cfn)
      (cfn) edge[->] (c);
\draw
   (addn) edge [->]  (addm)
   (addm) edge [->]  (covm)
   (covm) edge [->]  (nonn)
   (nonn) edge [->]  (cfn);
\draw (addm) edge [->] (b)
      (b)  edge [->] (nonm);
\draw (covm) edge [->] (d)
      (d)  edge[->] (cfm);
\draw (b) edge [->] (d);

\draw[color=blue] (-1,5.6)--(5,5.6);
\draw[color=blue] (-0.2,2.5)--(-0.2,7.5);
\draw[color=blue] (2.5,2.5)--(2.5,7.5);
\draw[color=blue] (5,7.5)--(5,2.5);
\draw[color=blue] (7.2,7.5)--(7.2,2.5);
\draw[color=blue] (9.8,7.5)--(9.8,2.5);
\draw[color=blue] (5,4.4)--(11,4.4); 

\draw[circle, fill=yellow,color=yellow] (1.7,4) circle (0.4);
\draw[circle, fill=yellow,color=yellow] (3.1,4)   circle (0.4);
\draw[circle, fill=yellow,color=yellow] (1.8,6.3) circle (0.4);
\draw[circle, fill=yellow,color=yellow] (7.8,3.7) circle (0.4);
\draw[circle, fill=yellow,color=yellow] (5.6,6) circle (0.4);
\draw[circle, fill=yellow,color=yellow] (5.6,3.8) circle (0.4);
\draw[circle, fill=yellow,color=yellow] (3.1,6.3) circle (0.4);
\draw[circle, fill=yellow,color=yellow] (7.8,6) circle (0.4);
\draw[circle, fill=yellow,color=yellow] (10.5,6) circle (0.4);
\node at (10.5,6) {$\theta_9$};
\node at (7.8,6) {$\theta_8$};
\node at (5.6,3.8) {$\theta_5$};
\node at (3.1,6.3) {$\theta_4$};
\node at (1.7,4) {$\theta_1$};
\node at (3.1,4) {$\theta_3$};
\node at (1.8,6.3) {$\theta_2$};
\node at (7.8,3.7) {$\theta_7$};
\node at (5.6,6) {$\theta_6$};
}
\end{tikzpicture}
\caption{ Cicho\'n's maximum. The maximum number of values that Cicho\'n's diagram can assume.}
\label{cichonmax}
\end{figure}
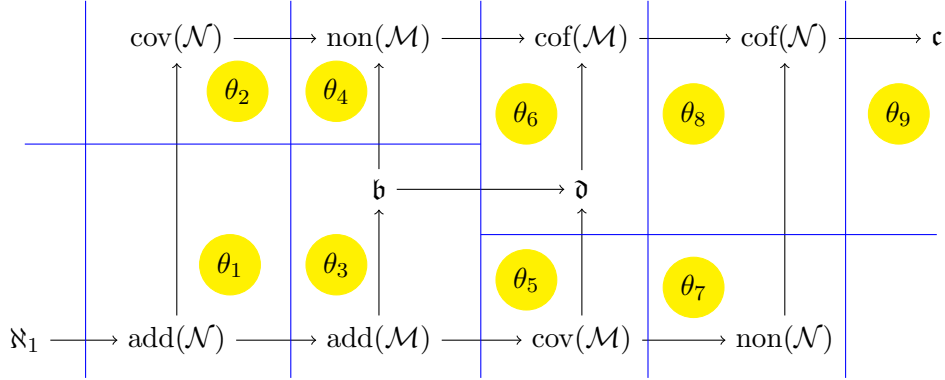

It turns out that the answer to~\autoref{Conj:Cichons} is positive and was proved by Goldstern, Kellner and Shelah~\cite{GKS}, who used four strongly compact cardinals to obtain the consistency of Cicho\'n's diagram divided into 10 different values, situation known as~\emph{Cicho\'n's maximum} (see~\autoref{cichonmax}). In this same direction. This was improved by Brendle, Mej\'ia, and the author~\cite{BCM} who produced 
a new method of matrix iterations with matrices of ultrafilters (called~\emph{${<}\mu$-uf-extendable matrix iterations}, see~\autoref{Defufextmatrix}), which extends the method from~\cite{GMS} to prove Cicho\'n's maximum modulo three strongly compact cardinals; finally, Goldstern, Kellner,  Mej\'ia, and Shelah~\cite{GKMS} proved that no large cardinals are needed for the consistency of Cicho\'n's maximum. 


Though we know how to separate all the entries of~\autoref{FigCichon}, the following is not known

\begin{question*}
Is it consistent that all the cardinals in~\autoref{FigCichon} (with the exception of the dependent values $\add(\Mwf)$ and $\cof(\Mwf)$) are pairwise different where $\cov(\Mwf) <\non(\Mwf)$?
\end{question*}

A natural question to ask about~~\autoref{diag:idealI} that motivate our interest in this work is the following.

\begin{question*}\label{Mproblem}
Is it consistent that all the four cardinals cardinal characteristics associated with $\Ewf$ in~\autoref{diag:idealI} are pairwise different?   
\end{question*}

In general, for any ideal $\Iwf$, there can be at most two instances of \autoref{Mproblem}, namely

\begin{itemize}
    \item[(A1)$_{\Iwf}$] $\add(\Iwf)<\cov(\Iwf)<\non(\Iwf)<\cof(\Iwf)$, and
    \mn
    \item[(A2)$_{\Iwf}$] $\add(\Iwf)<\non(\Iwf)<\cov(\Iwf)<\cof(\Iwf)$. 
\end{itemize}

The aim of this paper is to prove the following result:


\begin{teorema}{\emph{(\autoref{thm:4E2})}}\label{4E2}
Let $\theta_0\leq\theta\leq \mu\leq \nu$ be uncountable regular cardinals and let $\lambda$ be a cardinal such that $\nu\leq\lambda=\lambda^{{<}\theta}$. Then there is a ccc poset forcing~\autoref{SepofE}

\begin{figure}[ht!]
\centering
\begin{tikzpicture}
\small{
\node (aleph1) at (-1,3) {$\aleph_1$};
\node (addn) at (1,3){$\add(\Nwf)$};
\node (covn) at (1,7){$\cov(\Nwf)$};
\node (nonn) at (9,3) {$\non(\Nwf)$} ;
\node (cfn) at (9,7) {$\cof(\Nwf)$} ;
\node (addm) at (3.66,3) {$\add(\Mwf)$} ;
\node (covm) at (6.33,3) {$\cov(\Mwf)$} ;
\node (nonm) at (3.66,7) {$\non(\Mwf)$} ;
\node (cfm) at (6.33,7) {$\cof(\Mwf)$} ;
\node (b) at (3.66,5) {$\bfrak$};
\node (d) at (6.33,5) {$\dfrak$};
\node (c) at (11,7) {$\cfrak$};
\draw (aleph1) edge[->] (addn)
      (addn) edge[->] (covn)
      (covn) edge [->] (nonm)
      (nonm)edge [->] (cfm)
      (cfm)edge [->] (cfn)
      (cfn) edge[->] (c);
\draw
   (addn) edge [->]  (addm)
   (addm) edge [->]  (covm)
   (covm) edge [->]  (nonn)
   (nonn) edge [->]  (cfn);
\draw (addm) edge [->] (b)
      (b)  edge [->] (nonm);
\draw (covm) edge [->] (d)
      (d)  edge[->] (cfm);
\draw (b) edge [->] (d);

\draw[color=blue] (-1,5.6)--(5,5.6);
\draw[color=blue] (-0.2,2.5)--(-0.2,5.6);
\draw[color=blue] (2.5,2.5)--(2.5,5.6);
\draw[color=blue] (5,7.5)--(5,2.5);
\draw[color=blue] (5,4.4)--(11,4.4); 

\draw[circle, fill=yellow,color=yellow] (1.7,4.6) circle (0.4);
\draw[circle, fill=yellow,color=yellow] (3.1,4)   circle (0.4);
\draw[circle, fill=yellow,color=yellow] (2.5,6.3) circle (0.4);
\draw[circle, fill=yellow,color=yellow] (7.8,3.7) circle (0.4);
\draw[circle, fill=yellow,color=yellow] (7.8,5.5) circle (0.4);
\node at (1.7,4.6) {$\theta_0$};
\node at (3.1,4) {$\theta$};
\node at (2.5,6.3) {$\mu$};
\node at (7.8,3.7) {$\nu$};
\node at (7.8,5.5) {$\lambda$};
}
\end{tikzpicture}
\caption{Five values in Cicho\'n's diagram.}
\label{SepofE}
\end{figure}
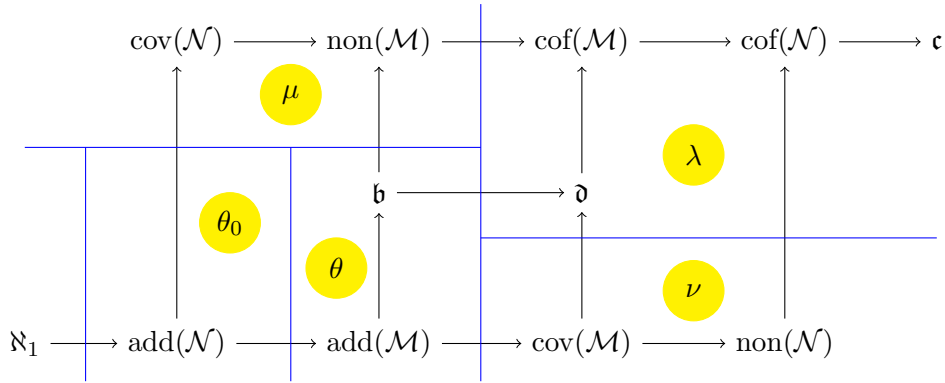
and $\add(\Ewf)=\theta$, $\non(\Ewf)=\mu$, $\cov(\Ewf)=\nu$, and $\cof(\Ewf)=\lambda$, i.e.\ \emph{(A2)$_{\Ewf}$.}
\end{teorema}
 
Concerning problems of this nature,  one of the earliest results is to due to Mej\'ia~\cite{M}, who used a forcing matrix iteration to get the consistency of $\onemainP_{\Nwf}$. Afterwards, classical preservation properties for matrix iterations were improved by the author joint with Mej\'ia \cite{CM} to force $\add(\Iwf_f)<\cov(\Iwf_f)<\non(\Iwf_f)<\cof(\Iwf_f)$ for any $f$ above some fixed $f^*$ (where $\Iwf_f$ is a~\textit{Yorioka ideal}, introduced by Yorioka~\cite{Yorioka} to show that no inequality between $\cof(\SNwf)$ and $\cfrak$ can be proven in ZFC), which solves $\onemainP_{\Iwf_f}$ but not for all $f$ at the same time. Later, $``\twomainP_\Mwf$ and $\onemainP_\Nwf"$ is consistent with ZFC\,$+$\,large cardinals is a consequence of Cicho\'n's maximum~\cite{GKS}. Afterwards, Brendle, the author and Mej\'ia~\cite{BCM} obtained that $\twomainP_\Mwf$ is consistent with ZFC (without large cardinals). 




    

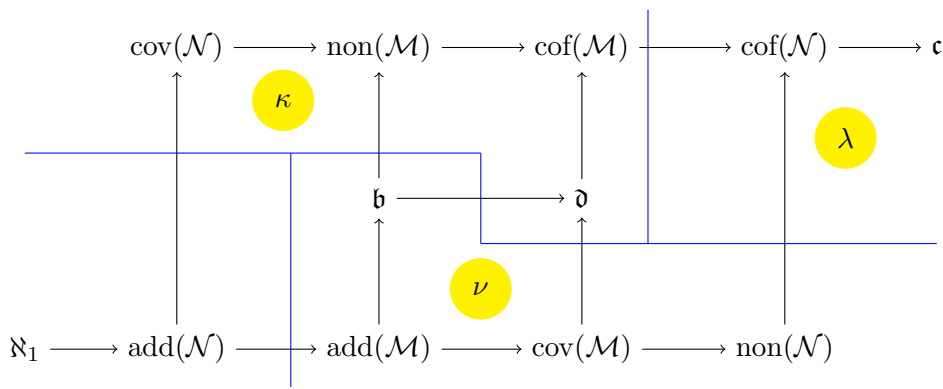
\begin{figure}[ht!]
\centering
\begin{tikzpicture}
\small{
\node (aleph1) at (-1,3) {$\aleph_1$};
\node (addn) at (1,3){$\add(\Nwf)$};
\node (covn) at (1,7){$\cov(\Nwf)$};
\node (nonn) at (9,3) {$\non(\Nwf)$} ;
\node (cfn) at (9,7) {$\cof(\Nwf)$} ;
\node (addm) at (3.66,3) {$\add(\Mwf)$} ;
\node (covm) at (6.33,3) {$\cov(\Mwf)$} ;
\node (nonm) at (3.66,7) {$\non(\Mwf)$} ;
\node (cfm) at (6.33,7) {$\cof(\Mwf)$} ;
\node (b) at (3.66,5) {$\bfrak$};
\node (d) at (6.33,5) {$\dfrak$};
\node (c) at (11,7) {$\cfrak$};
\draw (aleph1) edge[->] (addn)
      (addn) edge[->] (covn)
      (covn) edge [->] (nonm)
      (nonm)edge [->] (cfm)
      (cfm)edge [->] (cfn)
      (cfn) edge[->] (c);
\draw
   (addn) edge [->]  (addm)
   (addm) edge [->]  (covm)
   (covm) edge [->]  (nonn)
   (nonn) edge [->]  (cfn);
\draw (addm) edge [->] (b)
      (b)  edge [->] (nonm);
\draw (covm) edge [->] (d)
      (d)  edge[->] (cfm);
\draw (b) edge [->] (d);

\draw[color=blue] (-1,5.6)--(5,5.6);
\draw[color=blue] (2.5,2.5)--(2.5,5.6);
\draw[color=blue] (5,5.6)--(5,4.4);
\draw[color=blue] (7.2,7.5)--(7.2,4.4);
\draw[color=blue] (5,4.4)--(11,4.4); 

\draw[circle, fill=yellow,color=yellow] (2.4,6.3) circle (0.4);
\draw[circle, fill=yellow,color=yellow] (5,3.8) circle (0.4);
\draw[circle, fill=yellow,color=yellow] (9.8,5.8) circle (0.4);
\node at (9.8,5.8) {$\lambda$};
\node at (5,3.8) {$\nu$};
\node at (2.4,6.3) {$\kappa$};
}
\end{tikzpicture}
\caption{ The constellation of Cicho\'n's diagram forced in~\cite{Brshatt} where $\aleph_1<\nu<\kappa<\lambda$ with $\kappa$ and $\nu$ regular.}
    \label{Brendleshatared}
\end{figure}

On the other hand, Brendle has developed a new forcing iteration technique called \textit{Shattered iterations}~\cite{Brshatt} to construct a ccc poset forcing \autoref{Brendleshatared}. This result shows a possible way to separate Cichon's diagram into four parts with $\aleph_1<\cov(\Mwf)<\non(\Mwf)$. As a consequence, it was obtained the consistency of $\twomainP_{\Nwf}$. The consistency of $\onemainP_{\Mwf}$ is still open.

In another recent result in this direction, the author joint with Mej\'ia and \sloppy Rivera-Madrid \cite{CMR} prove the consistency of a weak version of~$\twomainP_{\SNwf}$, 
\[\add(\SNwf)=\non(\SNwf)<\cov(\SNwf)<\cof(\SNwf).\]
which is the first consistency result where more than two cardinal characteristics associated with $\SNwf$ are pairwise different, which is a weak version of  $\twomainP_{\SNwf}$. Later in \cite{cardona2020}, the author used a forcing matrix iteration to prove the consistency of a weak version of  $\onemainP_{\SNwf}$,
    \[\add(\SNwf)=\cov(\SNwf)<\non(\SNwf)<\cof(\SNwf).\]

A quite natural approach in order to prove~\autoref{4E2} is to use a matrix iteration of ccc forcings due to Blass and Shelah~\cite{B1S} who constructed the first example of a matrix iteration to prove the consistency of $\ufrak<\dfrak$ (where $\ufrak$ denote the~\emph{ultrafilter number}) and it was improved by Mej\'ia~\cite{M} to separate several cardinals on the left and right side of Cicho\'n's diagram simultaneously. The main problem with this approach is that it does not work to separate $\bfrak<\non(\Mwf)<\cov(\Mwf)<\dfrak$. 
To achieve this, we use the method \emph{ultrafilter-extendable} matrix iteration due to Brendle, the author, and Mej\'ia~\cite{BCM}, which extends the method from~\cite{GMS}.  Concrectely, we perfome a ${<}\theta$-uf-extendable matrix iteration  (see~\autoref{Defufextmatrix}) of height $\nu+1$ and length $\lambda\nu\mu$ (ordinal product)  of ccc partial orders, such that their iterands are determined by
\begin{enumerate}
    \item[(P1)] subforcings ccc of size ${<}\theta_0$,\smallskip
    \item[(P2)] subforcings of Hechler forcing of size ${<}\theta$, as well as through\smallskip
    \item[(P3)] vertically restricted random forcing to increase $\cov(\Nwf)$ along the iteration.
\end{enumerate}
This construction forces $\add(\Nwf)=\theta_0$, $\bfrak=\theta$, $\cov(\Nwf)=\non(\Mwf)=\mu$, $\cov(\Mwf)=\non(\Nwf)=\nu$, and $\dfrak=\lambda$. Therefore, $\add(\Ewf)=\theta$ and $\cof(\Ewf)=\lambda$ because $\add(\Ewf)=\add(\Mwf)$ and $\cof(\Ewf)=\cof(\Mwf)$. On the other hand, since $\Ewf\subseteq\Mwf\cap\Nwf$ (see~\cite[Lemma 2.6.1]{BJ}), $\non(\Ewf)\leq\min\{\non(\Mwf),\,\non(\Nwf)\}$ and $\max\{\cov(\Mwf),\,\cov(\Nwf)\}\leq\cov(\Ewf)$. Then $\non(\Ewf)\leq\mu$ and $\cov(\Ewf)\geq\nu$. 

How about $\mu\leq\non(\Ewf)$ and $\cov(\Ewf)\leq\nu$? To calculate $\mu\leq\non(\Ewf)$ and $\cov(\Ewf)\leq\nu$ we shall see (in~\autoref{uppbE}) that a lower bound of $\non(\Ewf)$ can be given in terms of localization cardinals, denoted by $\blc_{b,h}$ (see \autoref{def:localc}), as well as an upper bound of $\cov(\Ewf)$ is given in terms of the dual of $\blc_{b,h}$, denoted by $\dlc_{b,h}$. The key point is to iterate, in addition:
\begin{enumerate}
    \item[(P4)]  vertically restricted the forcing $\LOCor_{b,h}$ (see~\autoref{fnotions}) to increase $\blc_{b,h}$ along the iteration.
\end{enumerate}

Such restrictions work to guarantee that the matrix construction forces $\mu\leq\blc_{b,h}$ and $\dlc_{b,h}\leq\nu$. Hence $\non(\Ewf)\geq\mu$ and $\nu\leq\cov(\Ewf)$. 





 


\subsection*{Structure of the paper}\label{sec:structure}
We give a brief outline of this paper.

\noindent\autoref{preli}: We review the basic notation, relational systems and the Tukey order, our forcing notation, matrix iterations, and the classical preservation results for FS iterations.

\nonumber\autoref{charnoncovE}: We expand Bartoszy\'nski's and Shelah's characterizations about the covering and the uniformity of $\Ewf$.

\noindent\autoref{sec:filter-uflimits}: We review the notion of filter-linkedness from~\cite{mejiavert} (corresponding to $\Fr$-linkedness in~\autoref{Defufext}(1)) and its corresponding notions of linkedness and Knaster for poset from~\cite{BCM}, and we present some examples.


\noindent\autoref{apl}:  We prove that any complete Boolean algebra with a strictly-positive $\sigma$-additive measure preserves $\non(\Ewf)$ small. In particular, random forcing preserves $\non(\Ewf)$ small, which was proved by~Bartoszy\'nski and Shelah~\cite{BS1992}. We also prove~\autoref{WVcovnonE}, and ~\autoref{4E2}. 

\noindent\autoref{sec:disc}: We discuss some open problems.

\section{Preliminary facts}\label{preli}


\subsection{Notation} We introduce the following notation. The quantifiers $\forall^\infty n$ and $\exists^\infty n$ mean, respectively, ``for all but finitely many $n$" and ``there are infinitely many $n$". Denote by $\Lb$ the Lebesgue measure on $2^\omega$. For $f,g:\omega\to\mathrm{On}$, we write $f\leq^* g$ if $\forall^\infty n\, (f(n)\leq g(n))$, which is read \emph{$f$ is dominated by $g$}, and $f\neq^\infty g$ if\, $\forall^\infty n(f(n)\neq g(n))$, which is read \emph{$x$ is eventually different from $y$}. A \emph{slalom} is a function $\varphi:\omega\to[\omega]^{<\omega}$. For two functions $x$ and $\varphi$ with domain $\omega$, write $x\in^{*}\varphi$ for $\forall^{\infty}{n}(x(n)\in \varphi(n))$, which is read \textit{$\varphi$ localizes $x$}. For a function $b$ with domain $\omega$ and $h\in\omega^\omega$, denote $\seq_{<\omega}(b):=\bigcup_{n<\omega}\prod_{i<n}b(i)$, $\prod b:=\prod_{i<\omega}b(i)$ and $\Swf(b,h):=\prod_{i<\omega}[b(i)]^{\leq h(i)}$. For any set $A$, $\id_A$ denotes the identity function on $A$. Denote $\id:=\id_\omega$. A \textit{tree} is a $T\in\Pwf(2^{<\omega})$ such that $\forall s,t\in2^{<\omega}(t\in T\textrm{\ and } s\subseteq t\Rightarrow s\in T)$. For a tree $T$ denote $\lim(T):=\{x\in2^{\omega}:\forall n<\omega(x{\upharpoonright}n\in T)\}$. For $s\in2^{<\omega}$ denote $[s]:=\{x\in2^\omega: s\subseteq x\}$. Given $T\subseteq2^{n}$ for some $n<\omega$ define $[T]:=\bigcup_{t\in T}[s]$.


\begin{definition}\label{def:b_d}
 \begin{enumerate}[(i)]
    \item $\bfrak:=\min\{|F|:F\subseteq \omega^\omega\text{\ and }\neg\exists y\in \omega^\omega\, \forall x\in F\, (x\leq^*y)\}$, the \textit{unbounding number}.
    \mn
    \item $\dfrak:=\min\{|D|:D\subseteq \omega^\omega\text{\ and }\forall x\in \omega^\omega\, \exists y\in D\, (x\leq^* y)\}$, the \textit{dominating number}.
\end{enumerate}  
\end{definition}

We will be interested in the following cardinal characteristics, which were mentioned in the introduction:

\begin{definition}\label{def:localc}
Let $b=\la b(n):n<\omega\ra$ be a sequence of non-empty sets and let $h\in\omega^\omega$. Define the cardinals numbers $\blc_{b,h}$, $\dlc_{b,h}$, called \textit{localization cardinals}, as follows:
\begin{align*}
\cfa_{b,h} & :=\min\Big\{|D| :\, D\subseteq \Scal(b,h),\ \forall x \in \prod b\, \exists \varphi\in D\, (x\in^* \varphi)\Big\},\\[0.5ex]
\vfa_{b,h}&:=\min\Big\{|F| :\, F\subseteq \prod b,\ \neg\exists \varphi \in \Swf(b,h)\, \forall x \in F\, (x\in^* \varphi)\Big\},
\end{align*}
\end{definition}

The cardinals presented above have appeared in many contexts. For instance, these cardinals provide the well-known characterizations of the classical cardinal characteristics $\add(\Nwf)=\vfa_{\omega,h}$ and $\cof(\Nwf)=\cfa_{\omega,h}$ when $h$ diverges to infinity (here $\omega$ is interpreted as the constant sequence $\omega$), see~\cite[Thm. 2.3.9]{BJ}.



\subsection{Relational systems and Tukey order} In many cases, classical cardinal characteristics of the continuum, as well as those we deal with in this work, can be expressed by relational systems. 

\begin{definition}
A \textit{relational system} is a triple $\Rbf=\la X, Y, \sqsubset\ra$ where $X$ and $Y$ are non-empty sets and $\sqsubset$ is a relation contained in $X\times Y$. For $x\in X$ and $y\in Y$, $x\sqsubset y$ is often read $y\sqsubset$-\textit{dominates} $x$. 
\begin{enumerate}
    \item[(1)] A family $F\subseteq X$ is \textit{$\Rbf$-bounded} if there is a member of $Y$ that $\sqsubset$-dominates every member of $F$. Otherwise, we say that the set is \emph{$\Rbf$-unbounded}. 
    \mn
    \item[(2)] Dually, $D\subseteq Y$ is \textit{$\Rbf$-dominating} if every member of $X$ is $\sqsubset$-dominated by some member of $D$. 
\end{enumerate}
The following two cardinal characteristics are associated with $\Rbf$:
\begin{align*}
\bfrak(\Rbf)&:=\min\{|F|:F\subseteq X  \textrm{\ is\ }  \Rbf\textrm{-unbounded}\}\\
\dfrak(\Rbf)&:=\min\{|D|:D\subseteq Y \textrm{\ is\ } \Rbf\textrm{-dominating}\}.  
\end{align*} 
\end{definition}
These cardinals do not exist in general, for example, $\bfrak(\Rbf)$ does not exist iff $\dfrak(\Rbf)=1$; likewise, $\dfrak(\Rbf)$ does not exist iff $\bfrak(\Rbf)=1$.

Many classical cardinal characteristics can be expressed through relational systems.

\begin{example}\label{diset}
\begin{enumerate}
    \item[(1)] A \textit{preorder} is pair $\la S,\leq\ra$ where $S\neq\varnothing$ and $\leq$ is a relation on $S$ that satisfies reflexivity and transitivity. A \textit{directed preorder} is a preorder $\la S,\leq\ra$ that satisfies $\forall x, y\in S\,\exists z\in S\,(x\leq z\textrm{\ and\ }y\leq z)$.  As a relational system, $S=\la S, S\leq\ra$.
    \sn
    \item[(2)] For any ideal $\Iwf$ on $X$: 
    \sn
\begin{itemize}
    \item[(i)] $\Cbf_\Iwf:=\la X,\Iwf,\in\ra$, so $\bfrak(\Cbf_\Iwf)=\non(\Iwf)$ and $\dfrak(\Cbf_\Iwf)=\cov(\Iwf)$.
    \mn
    \item[(ii)] $\Iwf:=\la\Iwf,\subseteq\ra$ is directed, $\bfrak(\Iwf)=\add(\Iwf)$ and $\dfrak(\Iwf)=\cof(\Iwf)$.
\end{itemize}
\end{enumerate}
\end{example}

\begin{example}\label{Defantiloc}
Let $b$ be a function with domain $\omega$ such that $b(i)\neq\emptyset$ for all $i<\omega$, and let $h\in\omega^\omega$. Let $\Lc(b,h):=\la\prod b,\Scal(b,h),\in^*\ra$ be a relational system ($\Lc$ stands for \emph{localization}). Then  $\blc_{b,h}:=\bfrak(\Lc(b,h))$ and $\dlc_{b,h}:=\dfrak(\Lc(b,h))$.
\end{example}

In order to determine relations between cardinal characteristics,  the Tukey order is very useful.

\begin{definition}
Let $\Rbf':=\la X',Y',\sqsubset'\ra$ be another relational system. Say that \emph{$\Rbf$ is Tukey below $\Rbf'$}, denoted by $\Rbf\leqT\Rbf'$, if there are maps $F:X\to X'$ and $G:Y'\to Y$ such that, for any $x\in X$ and $y'\in Y'$, if $F(x)\sqsubset' y'$ then $x\sqsubset G(y')$. Say that \emph{$\Rbf$ and $\Rbf'$ are Tukey equivalent}, denoted by $\Rbf\eqT\Rbf'$, if $\Rbf\leqT\Rbf'$ and $\Rbf'\leqT\Rbf$. 

Note that $\Rbf\leqT\Rbf'$ implies $\bfrak(\Rbf')\leq\bfrak(\Rbf)$ and $\dfrak(\Rbf)\leq\dfrak(\Rbf')$.    
\end{definition}




\begin{example}\label{Charb-d}
Denote by $\Ior$ the set of interval partitions of $\omega$. Define the relational systems $\Dbf_1:=\la\Ior,\Ior,\sqsubseteq\ra$ and $\Dbf_2:=\la\Ior,\Ior,\ntriangleright\ra$ where, for any $I,J\in\Ior$,
   \[ I \sqsubseteq J \text{\ iff } \forall^\infty n\, \exists m\, (I_m\subseteq J_n); \quad
      I \ntriangleright J \text{\ iff } \forall^\infty n\, \forall m\, ( I_n\nsupseteq J_m).
   \]
   For each $I\in\Ior$ we define $f_I:\omega\to\omega$ and $I^{*2}\in\Ior$ such that $f_I(n):=\min I_{n}$ and $I^{*2}_n:=I_{2n}\cup I_{2n+1}$. For each increasing $f\in\omega^\omega$ define the increasing function $f^*:\omega\to\omega$ such that $f^*(0)=0$ and $f^*(n+1)=f(f^*(n)+1)$, and define $I^f\in\Ior$ such that $I_n^f:=[f^*(n),f^*(n+1))$.
   
   In Blass~\cite{blass} it is proved that $\Dbf\eqT\Dbf_1\eqT\Dbf_2$ where $\Dbf:=\la\omega^\omega,\omega^\omega,\leq^*\ra$. More concretely, for any increasing function $f, g\in\omega^\omega$ and $I\in\Ior$ defined by $I_n:=[f(n), f(n+1))$:
   \begin{enumerate}
       \item[(a)] $f=f_I\leq^* g$ implies $I\sqsubseteq I^g$. 
       \smallskip
       
       \item[(b)] $g\not\leq^* f=f_I$ implies $I^g \triangleright  I$.
   \end{enumerate}
\end{example}

\subsection{Forcing notions}\label{fnotions} We shall consider the following forcing notions. 


\subsection*{Random forcing} The \emph{random forcing} is defined by \[\Bor=\{T\subseteq2^{<\omega}:T\textrm{\ tree and }\Lb(\lim(T))>0\},\] ordered by $\subseteq$. Consider $\Pwf(2^{<\omega})$ as the Cantor space. Given a real $\delta\in(0,1)$ let \[\Bor_{\geq\delta}:=\{T\in\Bor:\Lb(\lim(T))\geq\delta\},\] which is a compact subspace in $\Pwf(2^{<\omega})$.
    
For $s\in2^{<\omega}$ and $m<\omega$ define 
\[\Bor(s,m):=\{T\in\Bor: \lim(T)\subseteq [s]\textrm{\ and\ } \Lb(\lim(T))\geq\Lb([s])(1-2^{-m-10}) \},\]
      which is also compact. Note that $\bigcup_{(s,m)\in2^{<\omega}\times\omega}\Bor(s,m)$ is dense in $\Bor$.

The next forcing notion was introduced by Brendle and Mej\'ia \cite{BrM} to increase $\blc_{b,h}$. 
 
\subsection*{Localization forcing} For $b,h\in\omega^\omega$ such that $\forall i<\omega(b(i)>0)$ and $h$ goes to infinity, the \emph{localization forcing} is defined by
\[\LOCor_{b,h}:=\{(p,n): p\in\Swf(b,h),\, n<\omega\textrm{\ and\ }\exists m<\omega\forall i<\omega(|p(i)|\leq m)\},\]
ordered by $(p',n')\leq(p,n)$ iff $n\leq n'$, $p'{\upharpoonright}n=p$, and $\forall i<\omega(p(i)\subseteq q(i))$. 

If $G\subseteq\LOCor_{b,h}$ is $\LOCor_{b,h}$-generic over $V$ and $\varphi_{\mathrm{gen}}:=\bigcup\{p:\exists n((n,p)\in G))\}$, then $\varphi_{\mathrm{gen}}\in\Swf(b, h)$ and  $x\in^*\varphi_{\mathrm{gen}}$ for every $x\in\prod b$ in the ground model.

For $s\in\seq_{<\omega}(b^+)$ where $b^+:=\la[b(n)]^{\leq h(n)}: n<\omega\ra$, and $m<\omega$, define
    \[\LOCor_{b,h}(s,m):=\{(p,n)\in\LOCor_{b, h}:\,s\subseteq p,\, n=|s|, \textrm{\ and\ }\forall i<\omega\,(|p(i)|\leq m)\}\]
    It is clear that $\LOCor_{b,h}=\bigcup_{(s,m)\in\seq_{<\omega}(b^+)\times\omega}\LOCor_{b,h}(s,m)$.

Though $\LOCor_{b,h}$ is not $\sigma$-centered it satisfies a property that lies between $\sigma$-centered and $\sigma$-linked, which was introduced by Kamo and Osuga~\cite{KO}:

\begin{definition}[{\cite{KO}}]\label{link}
Let $\rho,\pi \in \omega^\omega$. A forcing notion $\Pbb$ is \textit{$(\rho,\pi)$-linked} if there exists a sequence $\langle Q_{n,j}:n<\omega, j<\rho(n)\rangle$ of subsets of $\Pbb$ such that
\begin{enumerate}[(i)]
\item $Q_{n,j}$ is $\pi(n)$-linked for all $n<\omega$ and $j<\rho(n)$, and \smallskip
\item $\forall p\in \Pbb\forall^{\infty}{n<\omega}\exists{j<\rho(n)}(p\in Q_{n,j})$.
\end{enumerate}
Here, condition (ii) can be replaced by
\begin{enumerate}[(i')]
\setcounter{enumi}{1}
    \item $\forall p\in \Pbb\forall^{\infty}{n<\omega}\exists{j<\rho(n)}\exists q\leq p(q\in Q_{n,j})$
\end{enumerate}
because (i) and (ii') imply that the sequence of $Q'_{n,j}:=\{p\in\Pbb:\exists q\in Q_{n,j}(q\leq p)\}$ ($n<\omega$ and $j<\rho(n)$) satisfies (i) and (ii).
\end{definition}

\begin{lemma}[{\cite[Lemma 5.10]{BrM}}]\label{lem:rplinked}
Let $b,h,\pi,\rho \in\omega^\omega$ be non-decreasing functions with $b\geq 1$ and $h$ going to infinity. If $\{m_{k}\}_{k<\omega}$ is a non-decreasing sequence of natural numbers that goes to infinity and, for all but finitely many $k<\omega$, $k\cdot \pi(k)\leq h(m_k)$ and $k\cdot|[b(m_k-1)]^{\leq k}|^{m_k}\leq \rho(k)$, then $\LOCor_{b, h}$ is $(\rho,\pi)$-linked.
\end{lemma}

\begin{lemma}\label{centrlink}
If $\Pbb$ is $\sigma$-centered then $\Pbb$ is $(\rho,\pi)$-linked when $\rho:\omega\to\omega$ goes to $+\infty$.
\end{lemma}
\begin{PROOF}{\autoref{centrlink}}
See~\cite[Lemma 2.18]{CM}.
\end{PROOF}


\subsection*{Hechler forcing}
The \emph{Hechler forcing} is defined by   \[\Dor:=\{(s,f):\,s\in\omega^{<\omega}\textrm{\ and\ }f\in\omega^\omega\textrm{\ strictly increasing}\},\]
ordered by $(t,g)\leq(s,f)$ if $s\subseteq t$, $f\leq g$ and $f(i)\leq t(i)$ for all $i\in |t|\menos|s|$. Recall that $\Dor$ is $\sigma$-centered and it adds a real $d$ in $\omega^\omega$ which \textit{is dominating over the ground model reals in $\omega^\omega$}, which means that  $f\leq^* d$ for any $\forall{f}\in V\cap \omega^\omega$. Even more, $\Dor$ adds  Cohen reals.

The following is a generalization of the standard ccc poset that adds an eventually different real (see e.g.~\cite{KO,CM}).

\subsection*{ED forcing}  Fix $b:\omega\to\omega+1\menos\{0\}$ and $h\in\omega^\omega$ such that $\lim_{i\to+\infty}\frac{h(i)}{b(i)}=0$ (when $b(i)=\omega$, interpret $\frac{h(i)}{b(i)}$ as $0$).
      The \emph{$(b,h)$-ED (eventually different real) forcing} $\Eor^h_b$ is defined as the poset whose conditions are of the form $p=(s,\varphi)$ such that, for some $m:=m_{p}<\omega$, 
         \begin{enumerate}[(i)]
            \item $s\in\seq_{<\omega}(b)$, $\varphi\in\Swf(b,m\cdot h)$, and
            \mn
            \item $m\cdot h(i)<b(i)$ for every $i\geq|s|$,
         \end{enumerate}
      ordered by $(t,\psi)\leq(s,\varphi)$ iff $s\subseteq t$, $\forall i<\omega(\varphi(i)\subseteq\psi(i))$ and $t(i)\notin\varphi(i)$ for all $i\in|t|\menos|s|$.

      Put $\Eor^h_b(s,m):=\{(t,\varphi)\in\Eor^h_b:t=s\text{\ and }m_{(t,\varphi)}\leq m\}$ for $s\in\seq_{<\omega}(b)$ and $m<\omega$.

      Denote $\Eor_b:=\Eor^1_b$, $\Eor:=\Eor_\omega$, $\Eor_b(s,m):=\Eor^1_b(s,m)$, and $\Eor(s,m):=\Eor_\omega(s,m)$.
      
     It is not hard to see that $\Eor_b^h$ is $\sigma$-linked. Even more, whenever $b\geq^*\omega$, $\Eor_b^h$ is $\sigma$-centered. When $h\geq^*1$, $\Eor_b^h$ adds an eventually different real\footnote{More generally, $\Eor_b^h$ adds a real $e\in\prod b$ such that $\forall^\infty_{i<\omega}(e(i)\notin\varphi(i))$ for any $\varphi\in\Scal(b,h)$ in the ground model.} in $\prod b$.






\subsection{Preservation and examples} We review a few well-known results concerning the theory of preservation properties for FS iterations developed by Judah and Shelah~\cite{JS} and Brendle~\cite{Br}, which was generalized in~\cite[Sect.~4]{CM}.

\begin{definition}\label{DefRelSys}
Let $\Rbf=\la X,Y,\sqsubset\ra$ be a relational system and let $\theta$ be a cardinal.  
\begin{enumerate}[(1)]
  \item For a set $M$,
  \sn
  \begin{itemize}
      \item[(i)] An object $y\in Y$ is \textit{$\Rbf$-dominating over $M$} if $x\sqsubset y$ for all $x\in X\cap M$.
      \mn
      \item[(ii)] An object $x\in X$ is \textit{$\Rbf$-unbounded over $M$} if it $\Rbf^\perp$-dominating over $M$, that is, $x\not\sqsubset y$ for all $y\in Y\cap M$. 
  \end{itemize}
  \sn
  \item A family $D\subseteq Y$ is \emph{strongly $\theta$-$\Rbf$-dominating} if $|D|\geq\theta$ and, for any $x\in X$, $|\{y\in D : x\not\sqsubset y\}|<\theta$.
  \mn
  \item A family $F\subseteq X$ is \emph{strongly $\theta$-$\Rbf$-unbounded} if it strongly $\theta$-$\Rbf^\perp$-dominating, that is, $|F|\geq\theta$ and, for any $y\in Y$, $|\{x\in F : x\sqsubset y\}|<\theta$.
\end{enumerate}
\end{definition}

\begin{lemma}\label{StrUnb-effect}
   In the context of \autoref{DefRelSys}, assume that $F\subseteq X$ is a strongly $\theta$-$\Rbf$-unbounded family. Then:
   \begin{itemize}
       \item[(i)] $F$ is $\Rbf$-unbounded, in particular, $\bfrak(\Rbf)\leq|F|$.\footnote{The case $\theta\leq1$ is uninteresting and trivial. When $\theta=0$, the existence of a strongly $\theta$-$\Rbf$-unbounded family implies $Y=\emptyset$, so $\bfrak(\Rbf)=0$; when $F$ is a strongly $1$-$\Rbf$-unbounded family, $F\neq\emptyset$ and, for any $y\in Y$, no member of $F$ is $\sqsubset$-dominated by $y$, so $\bfrak(\Rbf)\leq1$.}
       \mn
       \item[(ii)] Whenever $\theta$ is regular, $|F|\leq\dfrak(\Rbf)$.
   \end{itemize}
\end{lemma}

The previous fact is the reason why strongly unbounded families are used to obtain upper bounds of $\bfrak(\Rbf)$ and lower bounds of $\dfrak(\Rbf)$, so their preservation in forcing extensions helps to force values for such cardinals.

The following two definitions are the central concepts for preservation of strongly unbounded families.

\begin{definition}
Say that $\Rbf=\langle X,Y,\sqsubset\rangle$ is a \textit{Polish relational system (Prs)} if the following is satisfied:
\begin{enumerate}[(i)]
\item $X$ is a perfect Polish space,
\mn
\item $Y$ is a non-empty analytic subspace of some Polish space $Z$ and
\mn
\item $\sqsubset\cap(X\times Z)=\bigcup_{n<\omega}\sqsubset_{n}$ where $\la\sqsubset_{n}\ra_{n<\omega}$  is some increasing sequence of closed subsets of $X\times Z$ such that, for any $n<\omega$ and for any $y\in Y$,
$(\sqsubset_{n})^{y}=\{x\in X:x\sqsubset_{n}y \}$ is closed nowhere dense.
\end{enumerate}

By (iii), $\la X,\Mwf(X),\in\ra$ is Tukey-Galois below $\Rbf$ where $\Mwf(X)$ denotes the $\sigma$-ideal of meager subsets of $X$. Therefore, $\bfrak(\Rbf)\leq \non(\Mwf)$ and $\cov(\Mwf)\leq\dfrak(\Rbf)$.
\end{definition}

\begin{definition}[Judah and Shelah {\cite{JS}}]\label{JS}
Let $\Rbf=\la X, Y, \sqsubset\ra$ be a Prs and let $\theta$ be a cardinal. A poset $\Por$ is \textit{$\theta$-$\Rbf$-good} if, for any $\Por$-name $\dot{h}$ for a member of $Y$, there is a non-empty $H\subseteq Y$ (in the ground model) of size ${<}\theta$ such that, for any $x\in X$, if $x$ is $\Rbf$-unbounded over  $H$ then $\Vdash x\not\not\sqsubset \dot{h}$.
\end{definition}

\begin{lemma}[{\cite[Lemma~4]{M}}]
If $\Rbf$ is a Prs and $\theta$ is an uncountable regular cardinal then any poset of size ${<}\theta$ is $\theta$-$\Rbf$-good. In particular, Cohen forcing is $\Rbf$-good.
\end{lemma}

\subsection{Examples of preservation properties}

Now, we present the instances of Prs and their corresponding good posets that we use in our applications.

\begin{enumerate}[(1)]\label{ExamPresPro}

    \item\label{Pres(null-cov)}  Define $\Omega_n:=\{a\in [2^{<\omega}]^{<\aleph_0}:\Lb(\bigcup_{s\in a}[s])\leq 2^{-n}\}$ (endowed with the discrete topology) and put $\Omega:=\prod_{n<\omega}\Omega_n$ with the product topology, which is a perfect Polish space. For every $x\in \Omega$ denote 
    \begin{itemize}
        \item $N_{x}:=\bigcap_{n<\omega}\bigcup_{s\in x(n)}[s]$, which is clearly a Borel null set in $2^{\omega}$.
        \item $N_{x}^{*}:=\bigcup_{n<\omega}\bigcap_{s\in x(n)}[s]$, which is clearly a $F_\sigma$ null set in $2^{\omega}$.
    \end{itemize}

    \begin{enumerate}
    \item Define the Prs $\Cn:=\la \Omega, 2^\omega, \sqsubset\ra$ where $x\sqsubset z$ iff $z\notin N_{x}$. Recall that any null set in $2^\omega$ is a subset of $N_{x}$ for some $x\in \Omega$, so $\Cn$ and $\Cbf_\Nwf^\perp$ are Tukey-Galois equivalent. Hence, $\bfrak(\Cn)=\cov(\Nwf)$ and $\dfrak(\Cn)=\non(\Nwf)$. 

     Any $\mu$-centered poset is $\mu^+$-$\Cn$-good (see e.g.~\cite{Br}). In particular, $\sigma$-centered posets are $\Cn$-good.\smallskip
    

      \item Define the Prs $\Ue:=\la 2^\omega,\Omega, \sqsubset^*\ra$ where $x\sqsubset^* z$ iff $x\in N_{z}^{*}$. Recall that any F$_{\sigma}$ null set in $2^\omega$ is a subset of $N_{x}^{*}$ for some $x\in \Omega$, so $\Ue$ and $\Cbf_\Ewf$ are Tukey-Galois equivalent. Hence, $\bfrak(\Ue)=\non(\Ewf)$ and $\dfrak(\Ue)=\cov(\Ewf)$. On the other hand, ramdon forcing is $\Ue$-good (see~\autoref{thm:pres_non(E)} in this work).
   \end{enumerate}
   \mn

    \item\label{Pres(non-mea)} Consider the Polish relational system $\Ed:=\la\omega^\omega,\omega^\omega,\neq^*\ra$. By~\cite[Thm.~2.4.1 \& Thm.~2.4.7]{BJ},  $\bfrak(\Ed)=\non(\Mwf)$ and $\dfrak(\Ed)=\cov(\Mwf)$.
    \mn
    \item\label{Pres(unb)} The relational system $\Dbf:=\la\omega^\omega,\omega^\omega,\leq^{*}\ra$ is Polish and $\bfrak(\Dbf)=\bfrak$ and $\dfrak(\Dbf)=\dfrak$.
     Any $\mu$-$\mathrm{Fr}$-linked poset (see~\autoref{Defufext}) is $\mu^+$-$\Dbf$-good (see \cite[Thm. 3.30]{mejiavert}).
    \mn
    \item\label{Pres(uni-null)}  For $\Hwf\subseteq\omega^\omega$ denote $\Lc^*(\Hwf):=\langle\omega^\omega,\Swf(\omega,\Hwf),\in^{*}\rangle$ where \[\Scal(\omega, \Hwf):=\{\varphi:\omega\to[\omega]^{<\aleph_0}:\exists{h\in\Hcal}\,\forall{i<\omega}(|\varphi(i)|\leq h(i))\}.\]
    Any $\mu$-centered poset is $\mu^+$-$\Lc^*(\Hwf)$-good (see~\cite{Br,JS}) so, in particular, $\sigma$-centered posets are $\Lc^*(\Hwf)$-good. 
     
     If $\Hwf$ is countable and non-empty then $\Lc^*(\Hwf)$ is a Prs because $\Swf(\omega,\Hwf)$ is $F_\sigma$ in $([\omega]^{<\aleph_0})^\omega$. In addition, if $\Hwf$ contains a function that goes to infinity then
$\bfrak(\Lc^*(\Hwf))=\add(\Nwf)$ and $\dfrak(\Lc^*(\Hwf))=\cof(\Nwf)$ (see~\cite[Thm.~2.3.9]{BJ}). Moreover, if for any $h\in\Hwf$ there is an $h^\prime\in\Hwf$ such that $\frac{h}{h^\prime}$ coverge to $0$, then any Boolean algebra with a strictly
positive finitely additive measure is $\Lc^*(\Hwf)$-good (\cite{Ka}). In particular, any
subalgebra of random forcing is $\Lc^*(\Hwf)$-good.

\end{enumerate}

A fundamental property of $(\rho,\pi)$-linked poset is the following:

\begin{lemma}[{\cite[Lemma 5.14]{BrM}}]\label{linkedpresadd(N)}
 For any $\pi,\rho,g_0\in\omega^\omega$ with $\pi$ and $g_0$ going to $+\infty$, there is a $\leq^*$-increasing sequence $\Hcal=\{g_n:n<\omega\}$ with $\frac{g_{n+1}}{g_n}\to\infty$ such that any $(\rho,\pi)$-linked poset is  $\Lc^*(\Hcal)$-good.
\end{lemma}

The following results indicate that strongly unbounded families can be added with Cohen reals, and the effect on $\bfrak(\Rbf)$ and $\dfrak(\Rbf)$ by a FS iteration of good posets.

\begin{lemma}\label{lem:strongCohen}
Let $\mu$ be a cardinal with uncountable cofinality, $\Rbf=\la X, Y, \sqsubset\ra$ a Prs and let $\la\mathbb{P}_{\alpha}\ra_{\alpha<\mu}$ be a $\lessdot$-increasing sequence of $\cf(\mu)$-cc posets such that  $\Por_\mu=\limdir_{\alpha<\mu}\Por_{\alpha}$. If $\Por_{\alpha+1}$ adds a Cohen real $\dot{c}_\alpha\in X$ over $V^{\Por_\alpha}$ for any $\alpha<\mu$, then $\Por_{\mu}$ forces that $\{\dot{c}_\alpha:\alpha<\mu\}$ is a strongly $\mu$-$\Rbf$-unbounded family of size $\mu$.
\end{lemma}
\begin{PROOF}{\autoref{lem:strongCohen}}
See e.g.~\cite[Lemma~4.15]{CM}.
\end{PROOF}

\begin{theorem}\label{sizeforbd}
Let $\theta$ be an uncountable regular cardinal, $\Rbf=\la X, Y, \sqsubset\ra$ a Prs, $\pi\geq\theta$ an ordinal, and let $\Por_{\pi}=\langle\Por_{\alpha},\Qnm_{\alpha}:\alpha<\pi\rangle$ be a FS iteration such that, for each $\alpha<\pi$, $\Qnm_{\alpha}$ is a $\Por_{\alpha}$-name of a  non-trivial $\theta$-$\Rbf$-good $\theta$-cc poset. Then, $\Por_{\pi}$ forces that $\bfrak(\Rbf)\leq\theta$ and  $|\pi|\leq\dfrak(\Rbf)$.
\end{theorem}
\begin{PROOF}{\autoref{sizeforbd}}
See e.g.~\cite[Thm.~4.15]{CM} or~\cite[Thm.~3.6]{GMS}.
\end{PROOF}

\subsection{Matrix iterations}
Let us review some usual facts about matrix iterations in the context of what we call \emph{simple matrix iterations}. In this type of matrix iterations only restricted generic reals are added, and preservation properties behave very nicely.

Fix transitive models $M\subseteq N$ of $\thzfc$ and a Prs $\Rbf=\langle X,Y,\sqsubset\rangle$ coded in $M$.

\begin{definition}\label{DefCompM}
Given two posets $\Por\in M$  and $\Qor$ (not necessarily in $M$), say that \textit{$\Por$ is a \textit{complete subposet of $\Qor$} with respect to $M$}, denoted by $\Por \lessdot_{M}\Qor$, if $\Por$ is a subposet of $\Qor$ and every maximal antichain in $\Por$ that belongs to $M$ is also a maximal antichain in $\Qor$.
\end{definition}

\begin{figure}[ht]
\begin{center}
  \includegraphics{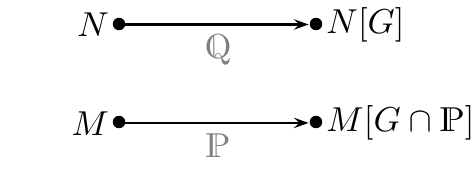}
  \caption{Generic extensions of pairs of posets ordered like $\Por\lessdot_M\Qor$.}
  \label{Figsinglestep}
\end{center}
\end{figure}

In this case, if $\Qor\in N$, then $\Por\lessdot_{M} \Qor$ implies that, whenever $G$ is $\Qor$-generic over $N$, $G \cap \Por$ is $\Por$-generic over $M$ and $M[G \cap \Por]\subseteq N[G]$ (see~\autoref{Figsinglestep}). When $\Por\in M$ it is clear that $\Por\lessdot_M\Por$.

The following results are related to preservation of $\Rbf$-unbounded reals along simple matrix iterations.

\begin{lemma}[{\cite[Lemma~11]{BrF}}, see also {\cite[Lemma~5.13]{mejia-temp}}]
   Assume that $\Por\in M$ is a poset. Then, in $N$,  $\Por$ forces that every $c\in X^N$ that is $\Rbf$-unbounded over $M$ is $\Rbf$-unbounded over $M^{\Por}$.
\end{lemma}

\begin{lemma}[{\cite{BrF}}]\label{parallellim}
   Assume that $\Por_{0,\pi}=\langle\Por_{0,\alpha},\Qnm_{0,\alpha}:\alpha<\pi\ra\in M$ and $\Por_{1,\pi}=\langle\Por_{1,\alpha},\Qnm_{1,\alpha}:\alpha<\pi\ra\in N$ are FS iterations such that, for any $\alpha<\pi$, if $\Por_{0,\alpha}\lessdot_M\Por_{1,\alpha}$ then $\Por_{1,\alpha}$ forces that $\Qnm_{0,\alpha}\lessdot_{M^{\Por_{0,\alpha}}} \Qnm_{1,\alpha}$. Then $\Por_{0,\alpha}\lessdot_M\Por_{1,\alpha}$ for any $\alpha\leq\pi$.

   In addition, if $\pi$ is limit, $c\in X^N$ and, for any $\alpha<\pi$, $\Por_{1,\alpha}$ forces (in $N$) that $c$ is $\Rbf$-unbounded over $M^{\Por_{0,\alpha}}$, then $\Por_{1,\pi}$ forces  that $c$ is $\Rbf$-unbounded over $M^{\Por_{0,\pi}}$.
\end{lemma}

\begin{definition}[Blass and Shelah {\cite{B1S}}]\label{Defmatit}
A \textit{matrix iteration} $\mbf$ consists of
\begin{enumerate}[(I)]
  \item a well order $I^\mbf$ and an ordinal $\pi^\mbf$, and
  
  \mn
  
  \item for each $i\in I^\mbf$, a FS iteration $\Pbb^\mbf_{i,\pi^\mbf}=\la\Pbb^\mbf_{i,\xi},\Qnm^\mbf_{i,\xi}:\xi<\pi^\mbf\ra$ such that, for any $i\leq j$ in $I^\mbf$ and $\xi<\pi^\mbf$, if $\Pbb^\mbf_{i,\xi}\lessdot\Pbb^\mbf_{j,\xi}$ then $\Pbb^\mbf_{j,\xi}$ forces $\Qnm^\mbf_{i,\xi}\lessdot_{V^{\Pbb^\mbf_{i,\xi}}}\Qnm^\mbf_{j,\xi}$.
\end{enumerate}
According to this notation, $\Pbb^\mbf_{i,0}$ is the trivial poset and $\Pbb^\mbf_{i,1}=\Qnm^\mbf_{i,0}$.
By \autoref{parallellim}, $\Pbb^\mbf_{i,\xi}$ is a complete subposet of $\Pbb^\mbf_{j,\xi}$ for all $i\leq j$ in $I^\mbf$ and $\xi\leq\pi^\mbf$.

We drop the upper index $\mbf$ when it is clear from the context. If $j\in I$ and $G$ is $\mathbb{P}_{j,\pi}$-generic over $V$ we denote  $V_{i,\xi}=V[G\cap\Pbb_{i,\xi}]$ for all $i\leq j$ in $I$ and $\xi\leq\pi$ . Clearly, $V_{i,\xi}\subseteq V_{j,\eta}$ for all $i\leq j$ in $I$ and $\xi\leq\eta\leq\pi$. The idea of such a construction is to obtain a matrix $\langle V_{i,\xi}:i\in I,\ \xi\leq\pi\rangle$ of generic extensions as illustrated in~\autoref{FigMatit}.
   
   
\begin{figure}[ht]
\begin{center}
  \includegraphics[width=\textwidth]{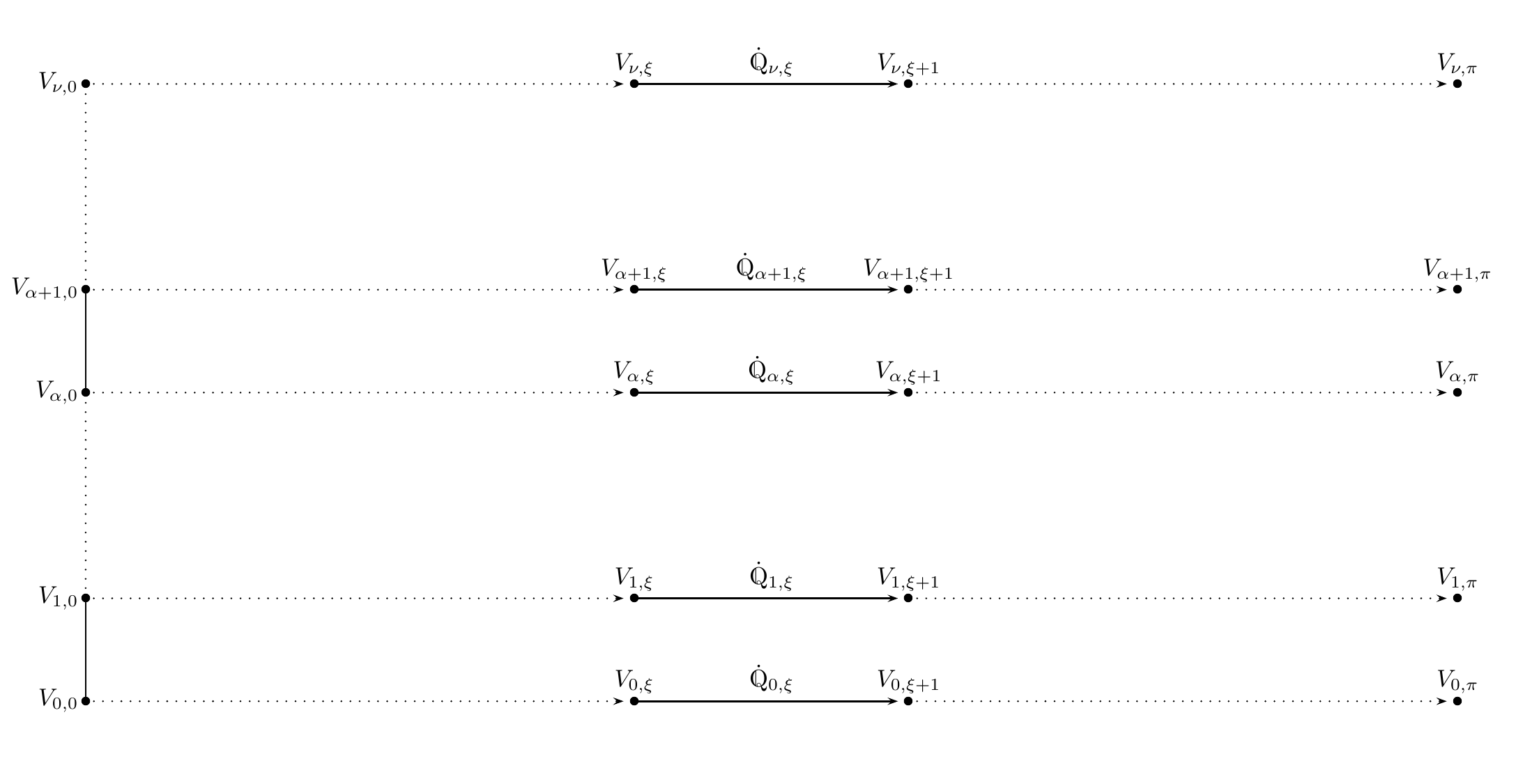}
  \caption{Matrix iteration with $I^\mbf=\nu+1$ where $\nu$ is an ordinal.}
  \label{FigMatit}
\end{center}
\end{figure}   
   
If $\xi\leq\pi$, $\mbf\frestr\xi$ (\emph{horizontal restriction}) denotes the matrix iteration with $I^{\mbf{\upharpoonright}\xi}=I$ and $\pi^{\mbf{\upharpoonright}\xi}=\xi$ where the FS iterations are the same as in (II) but restricted to $\xi$. On the other hand, for any $J\subseteq I$, $\mbf|J$ (\emph{vertical restriction}) denotes the matrix iteration with $I^{\mbf|J}=J$ and $\pi^{\mbf|J}=\pi$ where the FS iterations for $i\in J$ are exactly as in (II). 
\end{definition}

The following type of matrix iteration is the one we are going to deal with throughout the whole text.

\begin{definition}[Simple matrix iteration]\label{Defcohsys}
   A \emph{simple matrix iteration} $\mbf$ is a matrix iteration, composed additionally by a function $\Delta^\mbf:\pi^\mbf\to I^\mbf$, that satisfies: for each $\xi<\pi^\mbf$, there is a $\Por^\mbf_{\Delta^\mbf(\xi),\xi}$-name $\Qnm^\mbf_\xi$ of a poset such that, for each $i\in I^\mbf$,
         \[\Qnm^\mbf_{i,\xi}=\left\{\begin{array}{ll}
                 \Qnm^\mbf_\xi  & \text{if $i\geq\Delta^\mbf(\xi)$,}\\
                 \mathds{1} & \text{otherwise.}
                 \end{array}\right.\]
 
The upper index $\mbf$ is omitted when there is no risk of ambiguity.
\end{definition}
 
\begin{lemma}[{\cite{BrF}}, see also {\cite[Cor.~2.6]{mejiavert}}]\label{realint}
Let $\theta$ be an uncountable regular cardinal and let $\mbf$ be a simple matrix iteration as in~\autoref{Defcohsys} such that 
\begin{enumerate}
   \item[\emph{(i)}] $\gamma$, $I^{\mbf}=\gamma+1$, and $\cf(\gamma)\geq\theta$,\mn
   \item[\emph{(ii)}] $\gamma\notin\ran\Delta$, and\mn
   \item[\emph{(iii)}] for each $\xi<\pi$, $\Por_{\gamma,\xi}$ forces that $\Qnm_{\Delta(\xi),\xi}$ is $\theta$-cc.
\end{enumerate}
Then, for any $\xi\leq\pi$,
\begin{enumerate}
    \item[\emph{(a)}]  $\Por_{\gamma, \xi}$ is the direct limit of $\la\Por_{\alpha,\xi}:\alpha<\gamma\ra$, and
    \mn
    \item[\emph{(b)}] if $\dot{f}$ is a $\Por_{\gamma,\xi}$-name of a function from $\omega$ into $\bigcup_{\alpha<\gamma}V_{\alpha,\xi}$ then $\dot{f}$ is forced to be equal to a $\Por_{\alpha,\xi}$-name for some $\alpha<\gamma$.  In particular, the reals in $V_{\gamma,\xi}$ are precisely the reals in $\bigcup_{\alpha<\gamma}V_{\alpha,\xi}$.
\end{enumerate}  
\end{lemma}

\begin{theorem}[{\cite{BrF}}, see also {\cite[Thm.~10 \& Cor.~1]{M}}]\label{matsizebd}
   Let $\mbf$ be a simple matrix iteration with $I^\mbf=\gamma+1$ an ordinal, and let $\Rbf=\la X,Y,\sqsubset\ra$ be a Polish relational system coded in $V$. Assume that, for any $\alpha<\gamma$, there is some $\xi_\alpha<\pi$ such that $\Por_{\alpha+1,\xi_\alpha}$ adds a real $\dot{c}_\alpha\in X$ that is $\Rbf$-unbounded over $V_{\alpha,\xi_\alpha}$. 
   Then, for any $\alpha<\gamma$, $\Por_{\alpha+1,\pi}$ forces that $\dot{c}_{\alpha}$ is $\Rbf$-unbounded over $V_{\alpha,\pi}$.

   In addition, if $\Por_{\gamma,\pi}$ satisfies the hypothesis of~\autoref{realint} with $\gamma$ an ordinal of uncountable cofinality and $\theta=\cf(\gamma)$, and $f:\cf(\gamma)\to\gamma$ is increasing and cofinal, then  $\Por_{\gamma,\pi}$ forces that $\{\dot{c}_{f(\zeta)}:\zeta<\cf(\gamma)\}$ is a strongly $\cf(\gamma)$-$\Rbf$-unbounded family.
\end{theorem} 

 
 \section{Covering and uniformity of \texorpdfstring{$\Ewf$}{}}\label{charnoncovE}
This section is dedicated to expand Bartoszy\'nski's and Shelah's~\cite[Section 4]{BS1992} characterizations of the covering and the uniformity of $\Ewf$. What follows will prepare the way for the proof of these characterizations

The main combinatorial tool we are going to use is described below:

For $b\in\omega^\omega$ increasing let 
\[\Swf_b:=\bigg\{\la\varphi(n)\ra_{n<\omega}:\forall n<\omega\bigg(\varphi(n)\subseteq2^{[b(n),b(n+1))}\,\,\&\,\,\frac{|\varphi(n)|}{2^{b(n+1)-b(n)}}\leq\frac{1}{2^n}\bigg)\bigg\},\]
and 
\[\Swf_b^*:=\bigg\{\la\varphi(n)\ra_{n<\omega}:\forall n<\omega\bigg(\varphi(n)\subseteq2^{[b(n),b(n+1))}\bigg)\,\,\&\,\,\exists^\infty n\bigg(\frac{|\varphi(n)|}{2^{b(n+1)-b(n)}}\leq\frac{1}{2^n}\bigg)\bigg\}.\]

The motivation for this definition is the following result which will be useful. But we first begin remarking the following:

\begin{remark}\label{remarkE}
Let $b$ be a function with domain $\omega$ such that $b(i)\neq\emptyset$ for all $i<\omega$, and let $h\in\omega^\omega$. Notice that \[\limsup_{n\to\infty}\frac{h(n)}{|b(n)|}<1\textrm{\ iff\ }\forall A\in[\omega]^{\aleph_0}\left(\prod_{n\in A}\frac{h(n)}{|b(n)|}=0\right)\]
\end{remark}

\begin{lemma}\label{toolE}
For $b\in\omega^\omega$ and $\varphi\in\Swf_b^*$ such that $\limsup_{n\to\infty}\frac{|\varphi(n)|}{|b(n)|}<1$ the set
\[H_{b,\varphi}:=\{x\in2^\omega:\forall^\infty n\,(x{\upharpoonright}[b(n), b(n+1))\in\varphi(n))\}\]
belongs to $\Ewf$.
\end{lemma}
\begin{PROOF}{\autoref{toolE}}
Notice that $H_{b,\varphi}$ is a countable union of the closed null sets  
\[B_{b, \varphi}^m:=\{x\in2^\omega:\forall n\geq m (x{\upharpoonright}[b(n), b(n+1))\in\varphi(n))\}\textrm{\ for\ }m\in\omega.\] 
Indeed,  
\begin{align*}
\Lb(B_{b, \varphi}^m)=&\prod_{n\geq m}\Lb(\{x\in2^\omega:x{\upharpoonright}[b(n), b(n+1))\in\varphi(n)\})\\
=&\prod_{n\geq m}\frac{|\varphi(n)|}{2^{b(n+1)-b(n)}}\\
\leq&\prod_{n\in A}\frac{1}{2^n}=0 \textrm{\ where\ } A=\bigg\{n\geq m:\frac{|\varphi(n)|}{2^{b(n+1)-b(n)}}\leq\frac{1}{2^n}\bigg\},
\end{align*}
the latter inequality hold by~\autoref{remarkE}. Hence $\Lb(B_{b, \varphi}^m)=0$ for all $m<\omega$, so   
\[H_{b, \varphi}=\bigcup_{m<\omega} B_{b, \varphi}^m\] is an $F_\sigma$ null set and thus belongs to $\Ewf$.
\end{PROOF}

By the next lemma, we obtain a basis of $\Ewf$.

\begin{lemma}[{\cite[Thm. 4.3]{BS1992}}]\label{lem:combE}
Suppose that $C\in\Ewf$. Then there are $b\in\omega^\omega$ and $\varphi\in\Swf_b$ such that $C\subseteq H_{b, \varphi}$.
\end{lemma}
\begin{PROOF}{\autoref{lem:combE}}
Let us assume wlog that $C\subseteq2^\omega$ is a null set of type $F_\sigma$. Then $C$ can be written as $\bigcup_{n\in\omega}C_n$ where $\la C_n:n\in\omega\ra$ is an increasing family of closed sets of measure zero. 

Note that each $C_n$ is a compact set. It is easy to see that, if $K\subseteq2^{\omega}$ is a compact null set, then 
\[\forall\varepsilon>0\,\forall^{\infty}n\,\exists T\subseteq2^{n}\,(K\subseteq[T])\,\wedge\,\frac{|T|}{2^{n}}<\varepsilon).\]
Hence, we can define $b\in\omega^\omega$ increasing as follows: $b(0):=0$ and 
\[b(n+1):=\min\bigg\{m>b(n):\exists T_n\subseteq2^{m}\bigg(C_n\subseteq[T_n]\textrm{\ and\ }\frac{|T_n|}{2^m}\leq\frac{1}{4^{b(n)}}\bigg)\bigg\}.\]

Now define, for $n\in\omega$, \[\varphi(n):=\{s{\upharpoonright}[b(n), b(n+1)):s\in T_n\}.\]
It is clear that $|\varphi(n)|\leq |T_n|$. Hence, for every $n<\omega$,
\begin{equation*}\label{eqE}
  \frac{|\varphi(n)|}{2^{b(n+1)-b(n)}}\leq\frac{1}{2^{b(n)}}\leq\frac{1}{2^n}. 
\end{equation*}
Thus, $\la\varphi(n)\ra_{n<\omega}\in\Swf_b$, and by~\autoref{toolE} $H_{b, \varphi}\in\Ewf$. We also have $C\subseteq H_{b, \varphi}$. 
\end{PROOF}


\begin{lemma}\label{lem:mon}
Let $b_0,b_1\in\omega^\omega$ be increasing functions, let $b_1^*\in\omega^\omega$ be as in~\autoref{Charb-d}, and let $\varphi\in\Swf_{b_0}$.
\begin{itemize}
    \item[(1)] If $b_0\leq^*b_1$ then there is $\varphi_*\in\Swf_{b_1^*}$ such that $H_{b_0, \varphi}\subseteq H_{b_1^*, \varphi_*}$.
    \mn
    \item[(2)] If $b_1\not\leq^*b_0$ then there is  $\varphi^*\in\Swf_{b_1^*}^*$ such that $H_{b_0, \varphi}\subseteq H_{b_1^*, \varphi^*}$.
\end{itemize}\end{lemma}
\begin{PROOF}{\autoref{lem:mon}}
Let $I_n:=[b_0(n),b_0(n+1))$ and $I_n^*:=[b_1^*(n), b_1^*(n+1))$ for all $n<\omega$. 
\smallskip

\noindent(2): According to~\autoref{Charb-d}, if $b_1\not\leq^*b_0$, then
\[\exists^\infty n\, \exists m\, ( I_m\subseteq I_n^*).\]

This allows us to define $\varphi^*(n)$ for each $n\in\omega$ as follows: \[\varphi^*(n)
   := 
    \begin{cases}
    \big\{s\in 2^{I_n^*}:I_m\subseteq I_n^*\wedge s{\upharpoonright} I_m\in \varphi(m)\big\} & \textrm{if $m$ is minimal such that }\\
               & \textrm{$n\leq m$ and $I_m\subseteq I_n^*$;} \\ 
       2^{I_n^*} & \textrm{if such $m$ does not exist.}\\
    \end{cases}
\]
It remains to show that $\varphi^*\in\Swf_{b_1^*}^*$. To see this, it suffices to guarantee that there are infinitely many $n$ such that $m\geq n$, $I_m\subseteq I_n^*$.

Since $\exists^\infty m\,(b_0(m)<b_1(m))$, we may choose one $m$ such that $b_0(m+1)<b_1(m+1)$. Next, choose $n<\omega$ such that $b_0(m)\in I_n^*$, that is, $b_1^*(n)\leq b_0(m)<b_1^*(n+1)$. We have to distinguish two cases:

\noindent\textbf{Case 1:} $b_0(m+1)\leq b_1^*(n+1)$. Then $b_1(n)\leq b_1^*(n)\leq b_0(m)< b_0(m+1)\leq b_1^*(n+1)$, so $I_m\subseteq I_n^*$. On the other hand $b_1(n)<b_0(m+1)<b_1(m+1)$, so $n<m+1$, i.e., $n\leq m$.\smallskip

\noindent\textbf{Case 2:} $b_0(m+1)>b_1^*(n+1)$. We further distinguish two subcases:

\noindent\textbf{Subcase 2.1:} $b_0(m-1)\geq b_1^*(n)$. In this case, we have $I_{m-1}\subseteq I_n^*$. Then $b_1(n)\leq b_1^*(n+1)\leq b_0(m+1)<b_1(m+1)$, so $b_1(n)<b_1(m+1)$, which implies that $n<m+1$, that is, $n\leq m-1$. \smallskip

\noindent\textbf{Subcase 2.2:} $b_0(m-1)< b_1^*(n)$. We will obtain a contradiction by in this case.  Note that $m-1\leq b_0(m-1)<b_1^*(n)$, so $m\leq b_1^*(n)$. Then $b_1(m+1)\leq b_1(b_1^*(n)+1)=b_1^*(n+1)$.

On the other hand, $b_0(m+1)>b_1^*(n+1)$, so $b_1^*(n+1)<b_0(m+1)<b_1(m+1)\leq b_1^*(n+1)$. Hence, $b_1^*(n+1)<b_1^*(n+1)$, which is a contradiction.

We have proved that that there are infinitely many $n$, such that for some $m\geq n$, $I_m\subseteq I_n^*$. Hence $\varphi^*\in\Swf_{b_1^*}^*$. Indeed, for such $m$ and $n$, then
\[\frac{|\varphi^*(n)|}{2^{|I_n^*|}}=\frac{|\varphi(m)|}{2^{|I_m|}}\leq\frac{1}{2^m}\leq\frac{1}{2^n}.\]
Also, $H_{b, \varphi}\subseteq H_{b_1^*, \varphi^*}$.\smallskip

\noindent(1):  Once again, according to~\autoref{Charb-d}, note that if $b_0\leq^*b_1$ then
\[\forall^\infty n\, \exists m\, (I_m\subseteq I_n^*).\]

Now, for each $n\in\omega$ set $\varphi_*(n)$ by \[\varphi_*(n)
   := 
    \begin{cases}
    \big\{s\in 2^{I_n^*}:I_m\subseteq I_n^*\wedge s{\upharpoonright} I_m\in \varphi(m)\big\} & \textrm{if $m$ is minimal such that }\\
               & \textrm{$n\leq m$ and $I_m\subseteq I_n^*$;} \\ 
       \varnothing & \textrm{if such $m$ does not exist.}\\
    \end{cases}
\]
In a way similar to~\textbf{Case 1} in (1), it can be guaranteed that for all but finitely many $n$ there is some $m\geq n$ such that  $I_m\subseteq I_n^*$. This implies that $\varphi_*\in\Swf_{b_1^*}$  
and $H_{b, \varphi}\subseteq H_{b_1^*, \varphi_*}$.
\end{PROOF}

We use a variation of the relational system $\Lc(b,h)$ to prove a characterization of $\cov(\Ewf)$ in terms of the $\dfrak$-localization cardinals.

\begin{definition}
Let $b$ be a function with domain $\omega$ such that $b(i)\neq\emptyset$ for all $i<\omega$, and let $h\in\omega^\omega$. Define \[\Swf_*(b,h)=\left\{\varphi\in\prod_{n<\omega}\Pwf(b(n)):\forall n\,(\varphi(n)\subseteq b(n))\,\&\,\exists^\infty n\,(|\varphi(n)|\leq h(n))\right\}.\] 
Let $\loc_*(b,h):=\la\prod b,\Scal_*(b,h),\in^*\ra$ be a relational system. Denote  $\bloc_{b,h}:=\bfrak(\loc(b,h))$ and $\dloc_{b,h}:=\dfrak(\Lc(b,h))$. 
\end{definition}

For combinatorial purposes, we use the notion of closed null set in $\prod b$.

Denote $\Ewf(\prod b):=\{X\subseteq\prod b :X\text{\ is a closed measure zero}\}$. Likewise, we use the notation $\Ewf(\R)$ and $\Ewf([0,1])$.

When $b\nleq^* 1$, the map $F_b:\prod b\to[0,1]$ defined by
\[F_b(x):=\sum_{n<\omega}\frac{x(n)}{b(n)}\]
is a continuous onto function, and preserve measure. Hence, this map preserves sets between $\Ewf(\prod b)$ and $\Ewf([0,1])$ via images and pre-images, therefore, the value of the cardinal characteristics associated with $\Ewf$ do not depend on the space $\prod b$, neither on $[0,1]$.

\begin{lemma}\label{uppbE}
With the notation from the previous definition, if $b\in\omega^\omega$ and $\limsup_{n\to\infty}\frac{h(n)}{b(n)}<1$, then $\Cbf_\Ewf\leqT\loc_*(b,h)\leqT\loc(b,h)$. In particular, $\cov(\Ewf)\leq\dloc_{b,h}\leq\dlc_{b,h}$ and $\blc_{b,h}\leq\bloc_{b,h}\leq\non(\Ewf)$. 
\end{lemma}
\begin{PROOF}{\autoref{uppbE}}
Since $\Cwf_{\Ewf(2^\omega)}\eqT\Cwf_{\Ewf([0,1])}\eqT\Cwf_{\Ewf(\prod b)}$ by the above remarks, we just prove that $\Cbf_{\Ewf(\prod b)}\leqT\loc_*(b,h)$. We find functions $\Psi_-:\prod b\to\prod b$ and $\Psi_+:\Swf_*(b,h)\to\Ewf(\prod b)$ such that for any $x\in\prod b$ and $\varphi\in\Swf_*(b, h)$, if $\Psi_-(x)\in^*\varphi$, then $x\in\Psi_+(\varphi)$.

For $\varphi\in \Swf_*(b,h)$ put \[\Psi_+(\varphi):=\left\{x\in\prod b:\forall^\infty n\,(x(n)\in \varphi(n))\right\}.\] 
We only have to check that $\Psi_+(\varphi)\in\Ewf(\prod b)$. Indeed, setting $B_{\varphi,m}:=\{x\in\prod b:\forall n\geq m \,(x(n)\in\varphi(n))\}$, we have
\begin{align*}
\Lb(B_{\varphi,m})&=\prod_{n\geq m}\Lb\left(\left\{x\in\prod b:x(n)\in\varphi(n)\right\}\right)\\
&=\prod_{n\geq m}\frac{|\varphi(n)|}{b(n)}\\
&\leq\prod_{n\in A}^\infty\frac{h(n)}{b(n)}=0 \textrm{\ where\ } A=\bigg\{n\geq m:\frac{|\varphi(n)|}{b(n)}\leq\frac{h(n)}{b(n)}\bigg\}.
\end{align*}
where the last equality follows by~\autoref{remarkE}. Hence, $\Lb(B_{\varphi,m})=0$ for all $m\in\omega.$ Thus $H_\varphi=\{x\in\prod b:\forall^\infty n\,(x(n)\in\varphi(n))\}=\bigcup_{m<\omega} B_{\varphi,m}$ is an $F_\sigma$ null set and belongs to $\Ewf(\prod b)$.

It is clear that if $\forall^\infty n\,(x(n)\in \varphi(n))$ then $x\in \Psi_+(\varphi)$.
\end{PROOF}

\begin{theorem}\label{thm:chcovE}
\begin{multline*}
 \min\left(\{\dfrak\}\cup\left\{\dloc_{b, h}:b, h\in\omega^\omega\mathrm{\ and\ }\limsup_{n\to\infty}\frac{h(n)}{b(n)}<1\right\}\right)\leq\cov(\Ewf)\leq\\
 \leq\min\bigg\{\dloc_{b, h}:b, h\in\omega^\omega\mathrm{\ and\ }\limsup_{n\to\infty}\frac{h(n)}{b(n)}<1\bigg\}.   
\end{multline*}
In particular, if  $\cov(\Ewf)<\dfrak$ then \[\cov(\Ewf)=\min\bigg\{\dloc_{b, h}:b, h\in\omega^\omega\textrm{\ and\ }\limsup_{n\to\infty}\frac{h(n)}{b(n)}<1\bigg\}.\]

\end{theorem}
\begin{PROOF}{\autoref{thm:chcovE}}
By \autoref{uppbE} we know that \[\cov(\Ewf)\leq\min\bigg\{\dloc_{b, h}:b, h\in\omega^\omega\textrm{\ and\ }\limsup_{n\to\infty}\frac{h(n)}{b(n)}<1\bigg\}.\] 

It remains to prove that \[\cov(\Ewf)\geq\min\left(\{\dfrak\}\cup\left\{\dloc_{b, h}:b, h\in\omega^\omega\textrm{\ and\ }\limsup_{n\to\infty}\frac{h(n)}{b(n)}<1\right\}\right).\] 
To see this, assume $\kappa<\min\big(\{\dfrak\}\cup\big\{\dloc_{b, h}:b, h\in\omega^\omega\textrm{\ and\ }\limsup_{n\to\infty}\frac{h(n)}{b(n)}<1\big\}\big)$ and $\Jwf\subseteq\Ewf$ of size $\kappa$ . We will prove that $\bigcup\Jwf\neq2^\omega$

For each $A\in\Jwf$, by \autoref{lem:combE} there are $b_A\in\omega^\omega$ and $\varphi_A\in\Swf_{b_A}$ such that $A\subseteq H_{b_A, \varphi_A}$. Since $|\Jwf|<\dfrak$, we can find an increasing function $b$ such that $b\not\leq^* b_A$ for $A\in\Jwf$. By \autoref{lem:mon}, we obtain a family $\{\varphi_{A}^*:A\in\Jwf\}\subseteq\Swf^*_{b^*}$ such that $A\subseteq H_{b^*, \varphi_A^*}$. Observe that actually, this family corresponds to a family of slaloms in $\Swf_*(b', h)$ of size $<\min\{\dloc_{b, h}:b, h\in\omega^\omega\textrm{\ and\ }\limsup_{n\to\infty}\frac{h(n)}{b(n)}<1\}$ with $b'(n):=2^{b^*(n+1)-b^*(n)}$ and $h(n):=\frac{b'(n)}{2^n}$ for all $n$.

Since $\{\varphi_{A}:A\in\Jwf\}$ corresponds to a family of slaloms in $\Swf_*(b',h)$ of size ${\leq}\kappa$,  we can find an $x\in\prod b'$ such that, for all $A\in\Jwf$, $\exists^\infty n(x(n)\notin\varphi_A(n))$. Let $y\in2^\omega$ be such that $x(n)=y{\upharpoonright}[b'(n), b'(n+1))$. It is clear that $y\not\in H_{b^*, \varphi_{A}}$ for all $A\in\Jwf$, so $y\notin A$. 
\end{PROOF}

In a similar way, we can obtain a characterization of $\non(\Ewf)$.

\begin{theorem}
\begin{multline*}
\sup\bigg\{\bloc_{b, h}:b, h\in\omega^\omega\mathrm{\ and\ }\limsup_{n\to\infty}\frac{h(n)}{b(n)}<1\bigg\}\leq\non(\Ewf)\leq\\
  \sup\left(\{\bfrak\}\cup\left\{\bloc_{b, h}:b, h\in\omega^\omega\mathrm{\ and\ }\limsup_{n\to\infty}\frac{h(n)}{b(n)}<1\right\}\right).   
\end{multline*}
\end{theorem}
In particular, if $\non(\Ewf)>\bfrak$ then \[\non(\Ewf)=\sup\bigg\{\bloc_{b,h}:b, h\in\omega^\omega\mathrm{\ and\ }\limsup_{n\to\infty}\frac{h(n)}{b(n)}<1\bigg\}.\]


\section{Filter-linkedness and uf-limits}\label{sec:filter-uflimits}

\subsection{Filter-linkedness} To make the paper complete and self-contained we present a review of some results from~\cite[Section 3]{BCM} and~\cite{mejiavert}.\smallskip 

 Given a poset $\Por$, the $\Por$-name $\dot{G}$ usually denotes the canonical name of the $\Por$-generic set. If $\bar{p}=\la p_n:n<\omega\ra$ is a sequence in $\Por$, denote by $\dot{W}_\Por(\bar{p})$ the $\Por$-name of $\{n<\omega:p_n\in\dot{G}\}$. When the forcing is understood from the context, we just write $\dot{W}(\bar{p})$.\smallskip

Denote by $\Fr:=\{x\subseteq\omega:|\omega\menos x|<\aleph_0\}$ the \emph{Frechet filter}. A filter $F$ on $\omega$ is \emph{free} if $\Fr\subseteq F$. A set $x\subseteq\omega$ is \emph{$F$-positive} if it intersects every member of $F$. Denote by $F^+$ the family of $F$-positive sets. Note that $x\in\Fr^+$ iff $x$ is an infinite subset of $\omega$.

\begin{definition}[{\cite[Def. 3.1]{BCM}}]\label{Defufext}
Let $\Por$ be a poset, $F$ a free filter on $\omega$ and let $\mu$ be an infinite cardinal. 
\begin{enumerate}[(1)]
    \item A set $Q\subseteq\Por$ is \emph{$F$-linked} if, for any sequence $\bar{p}=\la p_n:n<\omega\ra$ in $Q$, there exists a $q\in\Por$ that forces $\dot{W}(\bar{p})\in F^+$.
    \mn
    \item A set $Q\subseteq\Por$ is~\emph{ultrafilter-linked}, abbreviated \emph{uf-linked}, if $Q$ is $D$-linked for any non-principal ultrafilter $D$ on $\omega$.
    \mn
    \item $\Por$ is \emph{$\mu$-$F$-linked} if $\Por=\bigcup_{\alpha<\mu}P_\alpha$ for some sequence $\la P_\alpha:\alpha<\mu\ra$ of $F$-linked subsets of $\Por$.
    \mn
    
    \item $\Por$ is \emph{$\mu$-$F$-Knaster} if $\forall A\in[\Por]^{\mu}\exists Q\in[A]^{\mu}(Q$ is $F$-linked in $\Por)$. 
    \mn
    \item The notions \emph{$\mu$-uf-linked} and~\emph{$\mu$-uf-Knaster} are defined similarly (by replacing $F$-linked by $uf$-linked).
\end{enumerate}
\end{definition}

\begin{remark}\label{RemKnaster}
\begin{enumerate}[(1)]
    \item Let $\Por$ be a poset and $Q\subseteq\Por$.  Note that a sequence $\la p_n:n<\omega\ra$ in $Q$ witnesses that $Q$ is \underline{not} $\Fr$-linked iff the set $\{q\in\Por:\forall^\infty n<\omega(q\perp p_n)\}$ is dense. \smallskip
    
    \item Any poset of size $\leq\mu$ is $\mu$-uf-linked because any singleton is uf-linked. 
\end{enumerate}

\end{remark}

\begin{lemma}[{\cite[Lemma~5.5]{mejiavert}}]\label{quasiuf}
If $\Por$ is a ccc poset then $Q\subseteq\Por$ is uf-linked iff it is Fr-linked. 
\end{lemma}

We now define a notion of convergence for sequences $\la x_n:n<\omega\ra$ in a topological space, indexed by $\omega$, with respect to an arbitrary filter $D$ on $\omega$.

\begin{definition}\label{Deffamlim}
  Let $D$ be an ultrafilter on $\Pwf(\omega)$, $X$ a  topological space. If $\bar{x}=\la x_n:n<\omega\ra$  is a sequence on $X$ and $x\in X$, we say that \emph{$\bar{x}$ $D$-converges to $x$} if, for every open neighborhood $U$ of $x$, $\{n<\omega:x_n\in U\}\in D$. Here, we also say that \emph{$x$ is a $D$-limit of $\bar{x}$}.
\end{definition}

A $D$-limit of a sequence $\bar{x}$ does not necessarily exist, nor will it be unique. But under  appropriate  assumptions,  there  is  an  unique $D$-limit.

\begin{lemma}\label{lem:top}
Let $X$ be a compact topological space. Assume $\bar{x}=\la x_n:n<\omega\ra$ is a sequence in $X$ and $D$  is an ultrafilter on $\Pwf(\omega)$. Then $\bar x$ has an  ultrafilter limit.
\end{lemma}
\begin{PROOF}{\autoref{lem:top}}
For a proof, see~\cite[Lemma 3.6]{BCM}.
\end{PROOF}

Note that there is at most one $D$-limit in Hausdorff spaces. In this case, we denote by $\lim^D_n x_n$ the ultrafilter limit of $\bar{x}$. Existence can always be guaranteed from compactness.

Our next goal is to show that the posets $\Bor$, $\Eor_b^h$ and $\LOCor_b^h$ are $\sigma$-uf-linked (see~\autoref{fnotions} for these posets).

\begin{example}
Fix $b:\omega\to\omega+1\menos\{0\}$ and $h\in\omega^\omega$ such that $\lim_{i\to+\infty}\frac{h(i)}{b(i)}=0$. Consider $\Pwf(\omega)$ as the Cantor space (homeomorphic to $2^{\omega}$), which is compact. Note that, for any $h'\in\omega^\omega$, $\Swf(b,h')$ is a compact subspace of $\Pwf(\omega)^\omega$ with the product topology, so for any $m<\omega$, $\Swf(b,m\cdot h)$ is a compact space. Therefore, if $D$ is an ultrafilter on $\omega$, $s\in\seq_{<\omega}(b)$ and $\bar{p}=\la p_n:n<\omega\ra$, $p_n=(s,\varphi_n)$, is a sequence in $\Eor^h_b(s,m)$, then the sequence $\la\varphi_n:n<\omega\ra$ has its $D$-limit $\varphi$ in $\Swf(b,m\cdot h)$ for any ultrafilter $D$ on $\Pwf(\omega)$. In this case, we say that \emph{$\lim^D_n p_n:=(s,\varphi)$ is the $D$-limit of $\bar{p}$}. Note that, for any $k<\omega$, $k\in\varphi(i)$ iff $\{n<\omega: k\in\varphi_n(i)\}\in D$.  
\end{example}

The fact that  $\Eor^h_b$ is $\sigma$-uf-linked is a consequence of the following.

\begin{lemma}[{\cite[Lemma~3.8]{BCM}}]\label{EDuflim}
Let $D$ be a non-principal ultrafilter on $\Pwf(\omega)$. If $G$ is $\Eor^h_b$-generic over $V$ then, in $V[G]$, $D$ can be extended to an ultrafilter $D^*$ on $\Pwf(\omega)\cap V[G]$ such that, for any $(s,m)\in\seq_{<\omega}(b)\times\omega$ and any sequence $\bar{p}\in\Eor_b^h(s,m)\cap V$ that has its $D$-limit in $G$, $\dot{W}(\bar{p})[G]\in D^*$.
\end{lemma}

The following result indicates that random forcing is $\sigma$-uf-linked. To this end, we define the following: For a Boolean algebra $\Bor$, say that $\mu:\Bor\to[0,1]$ is a \textit{strictly positive finitely additive (s.p.f.a.) measure} if it fulfills:
\begin{itemize}
    \item $\mu(\mathbf{1}_{\Bor}) = 1$, 
    \mn
    \item $\mu(a\vee a') = \mu(a)+\mu(a')$ for all $a,a'\in\Bor$ such that $a\wedge a'=0_{\Bor}$, and 
    \mn
    \item $\mu(a) = 0$ iff $a = \mathbf{0}_{\Bor}$.
\end{itemize}
Note that any Boolean algebra with a s.p.f.a. measure is ccc.

\begin{lemma}[{\cite[Lemma~3.29]{mejiavert}}]\label{mejiavert}
   Any complete Boolean algebra that admits a strictly-positive $\sigma$-additive measure is $\sigma$-uf-linked. In particular, any random algebra is $\sigma$-uf-linked.
\end{lemma}
\begin{PROOF}{\autoref{mejiavert}}
By \autoref{quasiuf} it suffices to prove that any such algebra is $\sigma$-$\Fr$-linked. Let $\Bor$ be a complete Boolean algebra that admits a strictly positive $\sigma$-additive measure $\mu$. For $m<\omega$, define
\[B_m=\bigg\{a\in\Bor:\mu(a)\geq\frac{1}{m+1}\bigg\}.\]
It is clear that $\Bor=\bigcup_{m<\omega}B_m$.  To finish the proof, it is enough to show that $B_m$ is Fr-linked. Assume the
contrary, so by \autoref{RemKnaster}(1) there is a sequence $\la a_n:n<\omega\ra$ in $B_m$ and a maximal antichain $\la a'_n:n<\omega\ra$ in $\Bor$ (where each $\mu(a'_n)>0$) such that, for each $n<\omega$, $a'_n\wedge a_k=\mathbf{0}_{\Bor}$ for all but finitely many $k<\omega$. Next, construct an increasing function
$g:\omega\to\omega$ such that $a'_n\wedge a_k=\mathbf{0}_{\Bor}$ for all $k\geq g(n)$. Find $n^*<\omega$ such that the measure of $a^*:=\bigvee_{n<n^*}a'_n$ is strictly larger than $1-\frac{1}{m+1}$, which is possible because $\mathbf{1}_{\Bor}=\bigvee_{n<\omega}a'_n$. Hence $\mu(a^*\wedge a)>0$ for any $a\in B_m$, but this contradicts that
$\mu(a^*\wedge a_k)=0$ for all $k\geq  g(n^*-1)$, which finishes the proof.
\end{PROOF}

\begin{example} Let $b,h\in\omega^\omega$ such that $\forall i<\omega(b(i)>0)$ and $h$ going to infinity. Note that for $m<\omega$, the set $[\omega]^{\leq m}:=\{a\subseteq\omega:|a|\leq m\}$ is a compact subspace of $\Pwf(\omega)$, so it is closed under uf-limits. Therefore, if $D$ is an ultrafilter on $\omega$, $s\in\seq_{<\omega}(b^+)$, $m<\omega$, and $\bar p=\la (p_n,|s|):n<\omega\ra$ is a sequence in $\LOCor_{b,h}(s,m)$, then the sequence $\bar p$ has 
 $D$-limit, and it is in $\LOCor_{b, h}(s, m)$.
\end{example}

Finally, we prove that $\LOCor_{b,h}$ is $\sigma$-uf-linked.

\begin{lemma}[{\cite[Remark~3.31]{mejiavert}}]\label{mejiavertLOC}
Let $b,h\in\omega^\omega$ such that $\forall i<\omega(b(i)>0)$ and $h$ goes to infinity, let $D$ be a non-principal ultrafilter on $\omega$ and $(s,m)\in\seq_{<\omega}(b^+)\times \omega$. If $\bar p=\la (p_n, |s|):n<\omega\ra$ is a sequence in $\LOCor_{b,h}(s,m)$ then there is a $(q,n)\in\LOCor_{b,h}(s,m)$ such that if $a\in D$, then $(q,n)$ forces that $\dot W(\bar p)\cap a\neq\varnothing$ is infinite.
\end{lemma}
\begin{PROOF}{\autoref{mejiavertLOC}}
Assume that $\bar p=\la (p_n,|s|):n<\omega\ra$ is a sequence in $\LOCor_{b,h}(s,m)$. For each $i\geq|s|$, since $[b(i)]^{\leq m}$ is finite, we can find $q(i) \in [b(i)]^{\leq m}$ and an $a_i\in D$
such that $p_n(i)=q(i)$ for any $n\in a_i$. To get $(q,n)\in\LOCor_{b,h}(s, m)$ we define  $n:=|s|$ and $q(i):= s(i)$ for all $i<n$. It remains to show that $(q,n)$ forces $\dot W(\bar p)\cap a\neq\varnothing$ whenever $a\in D$. To see this, assume that $(r,n^*)\leq (q,n)$ and
$a\in D$. Since $(r,n^*)\in\LOCor_{b,h}$, we can find $k, m_0 < \omega$ such that $k\geq n^*\geq n=|s|$, $|r(i)| \leq m_0$ for any $i<\omega$, and for any $i \geq k$, $m_0 + m \leq h(i)$ . Choose some $n_0\in\bigcap_{i<k} a_i\cap a$ (put $a_i:=\omega$ for $i < |s|$). Note that $p_{n_0}(i)=q(i)$ for all $i<k$.

Now we define $q'(i)$ by 
\[q'(i)=\left\{\begin{array}{ll}
                r(i)  & \text{if $i < k$,}\\
                 r(i)\cup p_{n_0}(i) & \text{if $i\geq k$.}
                 \end{array}\right.\]
Since $|q'(i)| \leq m_0 + m \leq h(i)$ for all $i\geq k$, $(q',k)$ is a condition in $\LOCor_{b,h}$. Moreover, $(q', k)$ is a condition stronger than $r$ and $p_{n_0}$, so it forces $n_0\in \dot W(\bar p)\cap a$.
\end{PROOF}

We conclude this section with a general result about the preservation of strongly $\Dbf$-unbounded families by $\Fr$-Knaster posets. 

\begin{theorem}[{\cite[Thm. 3.12]{BCM}}]\label{FrKnasterpresunb}
   If $\kappa$ is an uncountable regular cardinal then any $\kappa$-$\Fr$-Knaster poset preserves all the strongly $\kappa$-$\Dbf$-unbounded families from the ground model.
\end{theorem}

\section{Application of the ultrafilter-extendable matrix iteration}\label{apl}

The purpose of this section is to show~\autoref{WVcovnonE} and~\autoref{4E2}. As described in the introduction, to prove our main result we use an ultrafilter-extendable matrix iteration, which is defined as follows:

\begin{definition}[{\cite[Def.~4.1]{BCM}}]\label{Defufextmatrix}
   Let $\theta$ be an uncountable cardinal. A \emph{${<}\theta$-ultrafilter-extendable matrix iteration} (abbreviated \emph{${<}\theta$-uf-extendable}) is a simple matrix iteration $\mbf$ such that, for each $\xi<\pi^\mbf$, $\Por^\mbf_{\Delta^{\mbf}(\xi),\xi}$ forces that $\Qnm^\mbf_\xi$ is a ${<}\theta$-uf-linked poset.

   As in \autoref{Defcohsys}, we omit the upper index $\mbf$ when understood.
\end{definition}

\begin{theorem}[{\cite[Thm.~4.3]{BCM}}]\label{mainpres}
Let $\theta\geq\aleph_1$ regular cardinal and let $\Por_{\gamma,\pi}$ be a ${<}\theta$-uf-extendable matrix iteration. Then $\Por_{\alpha,\pi}$ is $\theta$-uf-Knaster for any $\alpha\leq\gamma$. In particular, it preserves any strongly $\theta$-$\Dbf$-unbounded family.
\end{theorem}


The main point for proving~\autoref{WVcovnonE} is the following lemma. 

\begin{lemma}[{\cite[Thm.~ 5.3]{BS1992}}]\label{thm:pres_non(E)}
Any complete Boolean algebra with a strictly-positive $\sigma$-additive measure $\mu$ is $\Ue$-good.
\end{lemma}
\begin{PROOF}{\autoref{thm:pres_non(E)}}
Let $\dot y$ be a $\Bor$-name for a member of $\Omega$. Since $\Bor$ is $\omega^\omega$-bounding, by~\autoref{lem:combE} and~\autoref{lem:mon} we may assume that $\Bor$ forces $N^*_{\dot y}\subseteq H_{b, \dot\varphi}$ for some $b\in\omega^\omega$ in the ground model and a $\Bor$-name $\dot\varphi$ of a member of $\Swf_b$, that is, 
\[N^*_{\dot y}\subseteq H_{b,\dot \varphi}=\{x\in2^\omega:\forall^\infty n\,(x{\upharpoonright}[ b(n), b(n+1))\in\dot\varphi(n))\}.\]
For $s\in 2^{[b(n), b(n+1))}$ set $B_{n,s}:=\ldbrack s\in\dot\varphi(n)\rdbrack\in\Bor$. 
Now, for $n<\omega$, define $\varphi(n)$ by 
\[\varphi(n):=\bigg\{s\in 2^{[b(n), b(n+1))}:\mu(B_{n,s})\geq2^{-\lfloor\frac{n}{2}\rfloor}\bigg\}.\]
We claim that 
\[\frac{|\varphi(n)|}{2^{b(n+1)-b(n)}}\leq2^{-\lfloor\frac{n}{2}\rfloor}.\]  
Suppose that for some $n_0\in\omega$,  \[\frac{|\varphi(n)|}{2^{b(n_0+1)-b(n_0)}}>2^{-\lfloor\frac{n_0}{2}\rfloor}.\] 
For $S\subseteq\varphi(n)$ set 
\[i(S):=\max\left\{|X|:X\subseteq S,\,\mu\bigg(\bigwedge_{s\in X}B_{n,s}\bigg)>0\right\}.\]
By~\cite[Prop. 1]{Kelley}, 
\[2^{-\lfloor\frac{n}{2}\rfloor}\leq\inf\{\mu(B_{n,s}):s\in\varphi(n)\}\leq\inf\bigg\{\frac{i(S)}{|S|}:\varnothing\subsetneq S\subseteq\varphi(n)\bigg\},\]
in particular, \[2^{b(n_0+1)-b(n_0)-n_0}<\frac{|\varphi(n_0)|}{2^{\lfloor\frac{n_0}{2}\rfloor}}\leq i(\varphi(n_0)).\]
Choose $X\subseteq\varphi(n_0)$ such that $|X|>2^{b(n_0+1)-b(n_0)-n_0}$ and $\mu\big(\bigwedge_{s\in X}B_{n_0,s}\big)>0$. Hence, $\bigwedge_{s\in X}B_{n_0,s}\Vdash X\subseteq\dot\varphi(n_0)$, so $\bigwedge_{s\in X}B_{n_0,s}\Vdash |\dot\varphi(n_0)|>2^{b(n_0+1)-b(n_0)-n_0}$, which is a contradiction.

Thus
\[H_{b,\varphi}=\{x\in2^\omega:\forall^\infty n\,(x{\upharpoonright}[b(n), b(n+1))\in\varphi(n))\}\]
is a member of $\Ewf$ by~\autoref{toolE}. Find $y\in\Omega$ such that $H_{b,\varphi}\subseteq N^*_y$.

To end the proof, let us argue that, if $x\not\in N^*_y$, then $\Vdash_\Bor x\not\in N^*_{\dot y}$. Suppose that $x\not\in N^*_y$. Towards a contradiction assume that $p\Vdash_\Bor x\in N^*_{\dot y}$ for some $p\in\Bor$. Since $\Vdash N^*_y\subseteq H_{b, \dot\varphi}$, we can assume wlog that there is some $m\in\omega$ such that $p\Vdash\forall n\geq m\,(x{\upharpoonright}[b(n), b(n+1))\in\dot \varphi(n))$ and $\mu(p)>2^{-m}$.

On the other hand, since $x\not\in N^*_y$ and $H_{b, \varphi}\subseteq N_y^*$, we can find an $n\geq 2m$ such that $x{\upharpoonright}[b(n), b(n+1))\not\in\varphi(n)$. In particular, 
\[\mu(B_{n,\,x{\upharpoonright}[b(n), b(n+1)})<\frac{1}{2^m}.\]
Now define $q:=p\smallsetminus B_{n,\,x{\upharpoonright}[b(n), b(n+1))}$. Then $\mu(q)>0$ and $q\Vdash x{\upharpoonright}[b(n), b(n+1))\not\in\dot \varphi(n)$, which is a contradiction.
\end{PROOF}


\begin{theorem}\label{weakversionE}
Let $\theta\leq\nu$ be uncountable regular cardinals and let $\lambda$ be a cardinal such that $\nu\leq\lambda=\lambda^{{<}\theta}$. Then there is a ccc poset forcing
\[\non(\Ewf)=\bfrak=\theta\leq\cov(\Nwf)=\non(\Mwf)=\cov(\Mwf)=\non(\Nwf)=\nu\leq\cov(\Ewf)=\dfrak=\lambda.\]
\end{theorem}
\begin{PROOF}{\autoref{weakversionE}}
We shall perform a FS iteration $\Por=\la\Por_\xi,\Qnm_\xi:\,  \xi<\pi\ra$ of length $\pi:=\lambda\nu$ (ordinal product) as follows. For each $\rho<\nu$ denote $\eta_\rho:=\lambda\rho$. We define the iteration at each $\xi=\eta_\rho+\varepsilon$ for $\rho<\nu$ and $\varepsilon<\lambda$ as follows:
\[\Qnm_\xi:=\left\{\begin{array}{ll}
        \dot\Bor\,(\textrm{random forcing})   & \text{if $\varepsilon=0$,}\\
        \Dor^{\dot N_\xi} & \text{if $\varepsilon>0$,}
    \end{array}\right.\]
where $\dot N_{\eta_\rho}$ is a $\Por_{\xi}$-name of a transitive model of $\thzfc$ of size ${<}\theta$ when $\varepsilon>0$. 

Additionally, by a book-keeping argument, we make sure that all such models $N_\xi$ are constructed such that, for any $\rho<\nu$:
\begin{equation}\label{claimcov}
\parbox{0.8\textwidth}{if $A\in V_{\eta_\rho}$ is a subset of $\omega^\omega$ of size ${<}\theta$, then there is some $0<\varepsilon<\lambda$ such that  $A\subseteq N_{\eta_\rho+\varepsilon}$.}   \tag{$\clubsuit$} 
\end{equation}
We prove that $\Por$ is as required. Clearly, $\Por$ forces $\cfrak=|\pi|=\lambda$.\smallskip


Note that all its iterands are $\theta$-$\Ue$-good (by~\autoref{thm:pres_non(E)}) and $\theta$-$\Dbf$-good. Thereore, by~\autoref{sizeforbd}, $\Por$ forces $\bfrak,\,\non(\Ewf)\leq\theta$, and $\dfrak,\,\cov(\Ewf)\geq\lambda$. Even more, since $\Por$ forces that $\cfrak\leq\lambda$, along with the former it is forced that $\cov(\Ewf)=\dfrak=\lambda$.\smallskip

Next, we show that $\bfrak\geq \theta$. To see this, let $F\subseteq\omega^\omega\cap V_{\pi}$ be a family of size ${<}\theta$.  By~\autoref{realint}, there is $\rho<\nu$ such that $F\in V_{\eta_\rho}$. Therefore, there is some $0<\varepsilon<\lambda$ such that $F\subseteq N_{\eta_{\rho}+\varepsilon}$ by~(\ref{claimcov}). Then $\Qnm_\xi=\Dor_{}^{N_\xi}$ adds a real that dominates all the reals in $ F$ where $\xi=\eta_\rho+\varepsilon$.\smallskip


On the other hand, since the FS iteration that determines $\Por$ has lenght of cofinality $\nu$ and $\nu$-cofinally many full random reals are added by $\Bor$, $\Por$ forces $\non(\Nwf)\leq\nu\leq\cov(\Nwf)$. Even more, $\Por$ adds $\nu$-cofinally many Cohen reals that form a strongly $\nu$-$\Ed$-unbounded family of size $\nu$ by~\autoref{lem:strongCohen}. Therefore, $\Por$ forces $\non(\Mwf)=\bfrak(\Ed)\leq\nu$. Even more, by~\autoref{matsizebd}, $\Por$ forces $\cov(\Mwf)=\dfrak(\Ed)\geq \nu$. Hence, $\Por$ forces that $\cov(\Nwf)=\non(\Mwf)=\cov(\Mwf)=\non(\Nwf)=\nu$.\smallskip

Since $\Por$ forces that $\bfrak=\theta$ and by~\autoref{covnonE}, $\theta=\min\{\bfrak,\non(\Nwf)\}\leq\non(\Ewf)$, we get $\Por$ forces $\theta\leq\non(\Ewf)$. Therefore, $\Por$ forces that $\non(\Ewf)=\theta$. 
\end{PROOF}

As a consequence we get 

\begin{corollary}[{\cite[Thm.~ 5.5 and~5.6]{BS1992}}]
It is consistent with ZFC that $\non(\Ewf)<\min\{\non(\Nwf), \non(\Mwf)\}$ and $\cov(\Ewf)>\max\{\cov(\Nwf),\cov(\Mwf)\}$.
\end{corollary}







    
    
    

Now, we are ready to prove~\autoref{4E2}.

\begin{theorem}\label{thm:4E2}
Let $\theta_0\leq\theta\leq \mu\leq \nu$ be uncountable regular cardinals and let $\lambda$ be a cardinal such that $\nu\leq\lambda=\lambda^{{<}\theta}$. Then there are increasing functions $b, h\in\omega^\omega$ and a ccc poset forcing 
\begin{multline*}
    \add(\Nwf)=\theta_0\leq\bfrak=\theta\leq\cov(\Nwf)=\non(\Mwf)=\blc_{b,h}=\mu\leq\\ \leq\cov(\Mwf)=\non(\Nwf)=\dlc_{b,h}=\nu\leq\dfrak=\lambda=\cfrak.
\end{multline*}
As an immediately consequence, we obtain that $\add(\Ewf)=\theta$, $\non(\Ewf)=\mu$, $\cov(\Ewf)=\nu$, and $\cof(\Ewf)=\lambda$.
\end{theorem}
\begin{PROOF}{\autoref{thm:4E2}} Fix a bijection $g=(g_0, g_1,g_2):\lambda\to\{0, 1\}\times\nu\times\lambda$ and a function $t:\nu\mu\to\nu$ such that $t(\nu\delta+\alpha)=\alpha$ for each $\delta<\mu$ and $\alpha<\nu$. Denote $\eta_\rho:=\nu+\lambda\rho$ for each $\rho<\nu\mu$.\smallskip

The desired poset is $\Cor_\lambda\ast\Por$ where $\Por$ is contructed in $V_{0,0}:=V^{\Cor_\lambda}$ from a ${<}\theta$-uf-extendable matrix iteration $\mbf$. \smallskip

By working in $V_{0, 0}$, construct by recursion increasing functions $\rho, \pi, h, b$ in $\omega^\omega$ such that for all but finitelly many $k<\omega$
\begin{enumerate}
    \item[(a)] $k\cdot\pi(k)\leq h(k)$, \smallskip
    \item[(b)] $k\cdot\big|[b(k-1)]^{\leq k}\big|^k\leq\rho(k)$, and \smallskip
    \item[(c)] $b(k)>2h(k)$.
\end{enumerate}
By (a) and (b), and~\autoref{lem:rplinked}, we obtain that $\LOCor_{b,h}$ is $(\rho, \pi)$-linked. Therefore, by~\autoref{linkedpresadd(N)}, there is a $\leq^*$-increasing sequence $\Hcal=\la g_n:n<\omega\ra$ such that $\LOCor_{b}^{h}$ is $\Lc^*(\Hwf)$-good. Notice that~(c) implies that  $\limsup_{n\to\infty}\frac{h(n)}{b(n)}<1$, so by~\autoref{uppbE} $\cov(\Ewf)\leq\dlc_{b,h}$ and $\blc_{b,h}\leq\non(\Ewf)$. \smallskip

Work in $V_{0,0}$.  Let $I^\mbf:=\nu+1$, $\pi^\mbf:=\nu+\lambda\nu\mu$,

\begin{enumerate}[(I)]
    \item for any $\alpha<\nu$, $\Delta^\mbf(\alpha):=\alpha+1$ and $\Qnm_{\Delta(\alpha),\alpha}^\mbf=\Cor$,
\end{enumerate}

and define the matrix iteration in the intervals of the form $[\eta_\rho,\eta_{\rho+1})$ as follows. For $\alpha<\nu$ choose 
\begin{enumerate}
   \item[(0)] an enumeration $\{\Qnm_{0,\alpha,\zeta}^\rho:\zeta<\lambda\}$ of all the nice $\Por_{\alpha,\eta_\rho}$-names for all the posets which underlining set is a subset of $\theta_0$ of size ${<}\theta_0$ and $\Vdash_{\Por_{\nu,\lambda\rho}}$``$\Qnm_{0,\zeta}^\rho$ is ccc"; and\smallskip
    
    \item[(1)] an enumeration $\{\Qnm_{1, \alpha,\,\zeta}^\rho:\zeta<\lambda\}$ of all the nice $\Por_{\alpha,\eta_\rho}$-names for all the $\sigma$-centered subposets of Hechler forcing of size ${<}\theta$.
\end{enumerate}

For $\xi\in[\eta_\rho,\eta_{\rho+1})$, set 
\begin{enumerate}[(I)]
\setcounter{enumi}{1}
    \item $\Delta^\mbf(\xi)=t(\rho)+1$ and $\Qnm_{\lambda\rho}^\mbf=\LOCor_{b,h}^{V_{\Delta(\xi),\xi}}$ when $\xi=\eta_\rho$,
    
    \mn

    \item  $\Delta^\mbf(\xi)=t(\rho)+1$ and $\Qnm_{\lambda\rho}^\mbf=\Bor^{V_{\Delta(\xi),\xi}}$ when $\xi=\eta_\rho+1$,  and 
     \mn
    \item   $\Delta(\xi)=g_0(\varepsilon)+2$ and $\Qnm_{\xi}^{\mbf}=\Qnm^{\rho}_{g(\varepsilon)}$ when $\xi=\eta_\rho+2+\varepsilon$ for some $\varepsilon<\lambda$. 
\end{enumerate}
Note that, for each $\xi<\eta_\rho$ with $\rho<\nu\mu$, $\Por_{\Delta(\xi), \xi}$ forces that $\Qnm_{\xi}$ is  $\sigma$-uf-linked. The case $\xi=\eta_\rho$ it follows by~\autoref{mejiavert}; when $\xi=\eta_\rho+1$ follows by~\autoref{mejiavertLOC}; and when $\xi=\eta_\rho+1+\varepsilon$ for some $\varepsilon<\lambda$ it follows by~\autoref{RemKnaster}(2). Then, according to \autoref{Defufextmatrix}, the above settles the construction of $\mbf$ as a ${<}\theta$-uf-extendable matrix iteration. Set $\Por:=\Por_{\nu,\pi}$, which is ccc.

We shall prove that $\Por$ forces what we want. First of all, note that $\Por$ is a $\theta$-$\uf$-Knaster poset by \autoref{mainpres}. Second, for any regular cardinal $\kappa\in[\theta,\lambda]$, $\Por$ preserves the strongly $\kappa$-$\Dbf$-unbounded familiy of size $\kappa$ added previously by $\Cor_\kappa$. In particular, $\Por$ forces $\bfrak\leq\theta$ and  $\lambda\leq\dfrak$. \smallskip

Observe that $\Por$ can be obtained by the FS iteration $\la\Por_{\nu,\xi},\Qnm_{\nu,\xi}:\xi<\pi^\mbf\ra$ and that all its iterands are $\theta_0$-$\Lc^*(\Hwf)$-good. Therefore, by~\autoref{sizeforbd}, $\Por$ forces  $\add(\Nwf)\leq\theta_0$.\smallskip

With the same argument as in~\autoref{weakversionE}, $\Por$ forces that 
$\bfrak\geq \theta$ and $\mathrm{MA}_{<\theta_0}$ (which implies $\add(\Nwf)\geq\theta_0$). \smallskip


We now prove that $\Por$ forces $\mu\leq\blc_{b,h},\,\cov(\Nwf)$ and $\dlc_{b,h},\,\non(\Nwf)\leq\nu$. For each $\rho<\nu\mu$ denote by $\dot \varphi_\rho\in\Swf(b,h)\cap V_{t(\rho)+1,\eta_\rho+1}$ the generic slalom over $V_{t(\rho)+1,\eta_\rho}$ added by $\Qnm_{t(\rho)+1,\eta_\rho}$, and denote by $r_{\rho}\in V_{t(\rho)+1,\eta_\rho+2}\cap 2^\omega$ the random real added by $\Qnm_{t(\rho)+1,\eta_\rho+1}$ over $V_{t(\rho)+1,\eta_\rho+1}$. Hence $\mu\leq\blc_{b,h}$ is consequence of the following: 
\begin{equation}\label{claimblc}
\parbox{0.8\textwidth}{In $V_{\nu,\pi}$, each family of reals in $\prod b$ of size $<\mu$ is localized by some slalom $\varphi_{\rho}$.}   \tag{$\spadesuit_0$} 
\end{equation}
Indeed, let $\dot F$ be a $\Por$-name for a family of reals in $\prod b$. By \autoref{realint} there are $\alpha<\nu$, and $\rho_0<\nu\mu$ such that $\dot F$ is a $\Por_{\alpha,\eta_{\rho_0}}$-name. By the definition of $t$, find a $\rho\in[\rho_0,\nu\mu)$ such that $t(\rho)=\alpha$, so $\Por_{t(\rho)+1,\eta_{\rho+1}}$ forces  that $\dot\varphi_\rho$ localizes all the reals from $\dot F$.\smallskip
 
Since $\{\varphi_\rho\,:\,0<\rho<\nu\mu\}$ is a family of reals of size $\leq\nu$, by (\ref{claimblc}) and \autoref{realint} any real in $V_{\lambda,\lambda\kappa}\cap\omega^\omega $ is localizated by $\varphi_\rho$. Hence $\Por$ forces $\dlc_{b,h}\leq \nu$. By a similar argument, we can prove
\begin{equation}\label{claimcovN}
\parbox{0.8\textwidth}{In $V_{\nu,\pi}$, for each family of Borel null sets of size $<\mu$, there is some $r_\rho$ that is not in its union.}   \tag{$\spadesuit_1$} 
\end{equation}
Hence $\Por$ forces $\mu\leq\cov(\Nwf)$ and $\non(\Nwf)\leq\nu$. On the other hand, $\Por$ adds $\mu$-cofinally many Cohen reals that form a strongly $\mu$-$\Ed$-unbounded family of size $\mu$ by~\autoref{lem:strongCohen}. Therefore, $\Por$ forces $\non(\Mwf)=\bfrak(\Ed)\leq\mu$. 

By \autoref{matsizebd}, $\Por$ forces $\cov(\Mwf)=\dfrak(\Ed)\geq \nu$. Therefore, $\Por$ forces $\cov(\Nwf)=\non(\Mwf)=\balc_{b,h}=\mu$, and $\cov(\Mwf)=\non(\Nwf)=\dalc_{b,h}=\nu$ (recall that $\blc_{b, h}\leq\non(\Mwf)$ and $\cov(\Mwf)\leq\dlc_{b, h}$). In addition, $\Por$ forces:\medskip

\noindent\underline{$\add(\Ewf)=\theta$:} because $\add(\Ewf)=\bfrak$~(\cite{BS1992}); \smallskip

\noindent\underline{$\non(\Ewf)=\mu$:} By~\autoref{uppbE} $\blc_{b,h}\leq\non(\Ewf)$, so $\mu\leq\non(\Ewf)\leq\min\{\non(\Mwf),\,\non(\Nwf)\}=\mu$; \smallskip

\noindent\underline{$\cov(\Ewf)=\nu$:} By~\autoref{uppbE} $\cov(\Ewf)\leq\dlc_{b,h}$, so $\nu=\max\{\cov(\Mwf),\,\cov(\Nwf)\}\leq\cov(\Ewf)\leq\nu$; and \smallskip

\noindent\underline{$\cof(\Ewf)=\lambda$:} because $\cof(\Ewf)=\dfrak$~(\cite{BS1992}).
\end{PROOF}

\section{Open problems}\label{sec:disc}


We conclude our work with a selection of open problems concerning~\autoref{Mproblem} and~\autoref{thm:4E2}, this selection being guided mainly by the author's interest. We choose
problems that seem to demand further development of new iteration techniques, so that their solution might give new light on the interplay between forcing and cardinal characteristics.

Though in~\autoref{thm:4E2} we separated the four cardinal characteristics of $\Ewf$, the uniformity of $\Ewf$ appears equal to the covering of $\Nwf$ as well as the covering of $\Ewf$ appears equal to the uniformity of $\Nwf$. In view of this, it  would  be  interesting  to  improve  our construction to allow these cardinals to be pairwise different, so we ask if we could separate the left side of Cicho\'n's diagram along with the uniformity and the covering of $\Ewf$.

\begin{question}\label{Qanother7}
Are each one the following statements consistent with ZFC?
\begin{multline}
\tag{1}\aleph_1<\add(\Nwf)<\bfrak<\cov(\Nwf)<\non(\Ewf)<\non(\Mwf)<\cov(\Mwf)\\
<\cov(\Ewf)=\non(\Nwf)=\dfrak=\cfrak.
\end{multline}
\begin{multline}
\tag{2}\aleph_1<\add(\Nwf)<\bfrak<\non(\Ewf)<\cov(\Nwf)<\non(\Mwf)<\cov(\Mwf)\\
<\cov(\Ewf)=\non(\Nwf)=\dfrak=\cfrak.
\end{multline}
\begin{multline}
\tag{3}\aleph_1<\add(\Nwf)<\cov(\Nwf)<\bfrak<\non(\Ewf)<\non(\Mwf)<\cov(\Mwf)\\
<\cov(\Ewf)=\non(\Nwf)=\dfrak=\cfrak.  
\end{multline}
\end{question}

In~\cite{KST}, it was constructed FAMS (finitely additive measures) along a FS (finite support) iteration to force
\[\aleph_1<\add(\Nwf)<\bfrak<\cov(\Nwf)<\non(\Mwf)<\cov(\Mwf)=\cfrak.
\]
Using the analog of this technique for ultrafilters, in~\cite{GMS} it was constructed an ultrafilter along a FS iteration to force
\[\aleph_1<\add(\Nwf)<\cov(\Nwf)<\bfrak<\non(\Mwf)<\cov(\Mwf)=\cfrak,
\]
which was improved later in~\cite{BCM} with the ultrafilter-extendable matrix iteration method used in this work. The latter seems not to be enough to force desired values to $\cov(\Ewf)$ and $\non(\Ewf)$ without increasing $\bfrak$ too much. One natural approach to solve this problem would be using FAMS along a matrix iteration. The main problem with this approach is that we do not know how to preserve $\non(\Ewf)$ in this context. 

In~\autoref{weakversionE} it was proved the consistency of 
\[\non(\Ewf)<\cov(\Nwf)=\non(\Mwf)=\cov(\Mwf)=\non(\Nwf)<\cov(\Ewf).\] For this it was essencial to prove that random forcing preserves $\non(\Ewf)$ small~(\autoref{thm:pres_non(E)}), originally proved by Bartoszy\'nski and Shelah~\cite{BS1992} to obtain the consistency of $\non(\Ewf)<\min\{\non(\Mwf),\non(\Nwf)\}$.

The key to what the lemma did to prove that random forcing preserves $\non(\Ewf)$ was the characterization of closed measure zero sets \autoref{lem:combE} and~\autoref{lem:mon}(a). Since eventually different real forcing does not add dominating reals, we may ask whether~\autoref{lem:combE} and~\autoref{lem:mon} would help to solve:

\begin{question}
Does eventually different real forcing preserve $\non(\Ewf)$~small?
\end{question}

But this may be much harder because $\sigma$-uf-linked does not imply $\Ue$-good, for example $\LOCor_{b,h}$ is $\sigma$-uf-linked and increases $\non(\Ewf)$. 

One positive answer to~\autoref{Qanother7} along with the method of submodels of~\cite{GKMS} would help solving:

\begin{question}
Is  Cichon's maximum along with $\non(\Ewf)$ and $\cov(\Ewf)$ assuming distinct values consistent with ZFC?
\end{question}

 


Regarding~\autoref{Mproblem}, the following questions are still open.


\begin{question}\label{Q:Mproblem}
Are each one of the following statements consistent with ZFC?
\begin{itemize}
    \item[(1)] \emph{$\onemainP_\Ewf$}.
    \sn
    \item[(2)] \emph{$\onemainP_\Mwf$}.
     \sn
    \item[(3)] \emph{$\onemainP_\SNwf$}.
     \sn
    \item[(4)] \emph{$\twomainP_\SNwf$}.
     \sn
    \item[(5)] \emph{$\twomainP_{\Iwf_f}$} for any $f\in\omega^\omega$.
\end{itemize}
\end{question}

FS iterations of ccc forcings will not work to solve~\autoref{Q:Mproblem}, with the exception of~(3), because any such iteration forces $\non(\Mwf)\leq\cov(\Mwf)$, so alternative constructions are required. A creature forcing method based on the notion of decisiveness \cite{KS09,KS12}, developed in \cite{FGKS,CKM} could work to solve~\autoref{Q:Mproblem}, but this method is restricted to $\omega^\omega$-bounding forcings, that is, results in $\dfrak=\aleph_1$. So it won't work for~(2) neither for~(1) because $\dfrak=\aleph_1$ implies $\cov(\Nwf)=\cov(\Ewf)$ by~\autoref{covnonE}. It would be interesting to discover how to modify these constructions to allow unbounded reals while ensuring decisiveness. Another alternative is Brendle's method of shattered iterations.

\subsection*{Acknowledgments} 

The author thanks Dr.\ Diego Mej\'ia for his valuable comments while working on this paper as well as several stimulating discussions during the elaboration of this paper. His help greatly appreciated.


{\small
\bibliography{left}
\bibliographystyle{alpha}
}

\end{document}